\documentclass{gtart}

\def\ifplaintex{\expandafter\ifx\csname documentclass\endcsname\relax}

\ifplaintex 
\else
\fi

\expandafter\ifx\csname beginpicture\endcsname\relax
\expandafter\ifx\csname documentclass\endcsname\relax
\input pictex \else
\input prepictex \input pictex \input postpictex \fi\fi

\def\gt{{\mathsurround=0pt\it $\cal G\mskip-2mu$eometry \&\ 
$\cal T\!\!$opology}}        

\def\gtp{{\mathsurround=0pt\it $\cal G\mskip-2mu$eometry \&\ 
$\cal T\!\!$opology $\cal P\!$ublications}}  

\def\lognumber#1{\def\thelognumber{#1}}
\def\volumenumber#1{\def\thevolumenumber{#1}}
\def\papernumber#1{\def\thepapernumber{#1}}
\def\volumeyear#1{\def\thevolumeyear{#1}}

\def\pagenumbers#1#2{\def\startpage{#1}\def\finishpage{#2}}
\def\published#1{\def\publishdate{#1}}
\def\proposed#1{\def\theproposer{#1}}
\def\seconded#1{\def\theseconders{#1}}
\def\received#1{\def\receiveddate{#1}}
\def\revised#1{\def\reviseddate{#1}}
\def\accepted#1{\def\accepteddate{#1}}

\long\def\asciiabstract#1{\long\def\theasciiabstract{#1}}

\def\shorttitle#1{\def\theshorttitle{#1}}

\let\\\par\let\thelognumber\relax
\let\thevolumenumber\relax\let\thepapernumber\relax
\let\thevolumeyear\relax\let\thesamplenumber\relax\let\startpage\relax
\let\finishpage\relax\let\publishdate\relax\let\receiveddate\relax
\let\reviseddate\relax\let\accepteddate\relax\let\theasciititle\relax
\let\theasciiauthors\relax
\let\theasciiabstract\relax
\let\theasciiemail\relax\let\theshortauthors\relax\let\theshorttitle\relax

\long\def\maketitlep{   

\count0=\startpage

\gt\hfill      
\beginpicture
\setcoordinatesystem units <0.33truein, 0.33truein> point at 2.2 0.9
\setplotsymbol ({$\cal G$})
\plotsymbolspacing=9truept
\circulararc 315 degrees from 0 1 center at 0 0
\setplotsymbol ({$\cal T$})
\circulararc 315 degrees from 1 -1 center at 1 0
\endpicture
\break
{\small\ifx\thesamplenumber\relax 
Volume \else Sample
\fi\thevolumenumber\ (\thevolumeyear)
\startpage--\finishpage\nl
Published: \publishdate}
\vglue 0.5truein plus 0.4fil minus 0.1truein

{\parskip=0pt\leftskip 0pt plus 1fil\def\\{\par\smallskip}{\ifplaintex\large
\else\Large\fi\bf\thetitle}\par\medskip}   

\vglue 0pt plus 0.1fil 

{\parskip=0pt\leftskip 0pt plus 1fil\def\\{\par}{\sc\theauthors}
\par\medskip}

\vglue 0pt plus 0.1fil 

{\small\parskip=0pt\let\newline\\
{\leftskip 0pt plus 1fil\def\\{\par}{\sl\theaddress}\par}
\expandafter\ifx\theemail\relax    
\relax\else\vglue 5pt plus 0.02fil minus 2pt\def\\{\stdspace{\rm 
and}\stdspace} 
\cl{Email:\stdspace\tt\theemail}\fi
\ifx\theurl\relax                  
\relax\else\vglue 5pt plus 0.02fil minus 2pt\def\\{\stdspace{\rm 
and}\stdspace}
\cl{URL:\stdspace\tt\theurl}\fi\par}

\vglue 7pt plus 0.3fil minus 3pt

{\bf Abstract}
\vglue 5pt plus 0.1fil minus 2pt

\theabstract

\vglue 7pt plus 0.3fil minus 3pt

{\bf AMS Classification numbers}\quad Primary:\quad \theprimaryclass

Secondary:\quad \thesecondaryclass

\vglue 5pt plus 0.3fil minus 2pt

{\bf Keywords}\quad \thekeywords

\vglue 10pt plus 0.5fil minus 5pt

{\small  Proposed: \theproposer\hfill Received: \receiveddate\nl
Seconded: \theseconders\hfill 
\ifx\reviseddate\relax                         
Accepted: \accepteddate                        
\else
Revised: \reviseddate                          
\fi}
\eject
}       

\let\maketitlepage\maketitlep
\let\maketitle\maketitlepage

\font\phead=cmsl9 scaled 950
\font=cmsl9 scaled 1050
\font\pnum=cmbx10 scaled 913
\font\lnum=cmbx10 
\font\pfoot=cmsl9 scaled 950
\font=cmsl9 scaled 1050
\ifplaintex
\headline{\vbox to 0pt{\vskip -4.5mm\line{\small\phead\ifnum
\count0=\startpage ISSN 1364-0380 (on line)
1465-3060 (printed) \hfill {\pnum\folio}\else\ifodd\count0\def\\{ }%
\ifx\theshorttitle\relax\thetitle\else\theshorttitle\fi\hfill{\pnum\folio}
\else\def\\{ and }{\pnum\folio}\hfill\ifx\theshortauthors\relax\theauthors
\else\theshortauthors\fi\fi\fi}\vss}}
\footline{\vbox to 0pt{\vglue 0mm\line{\small\pfoot\ifnum\count0=\startpage
\copyright\ \gtp\hfill\else
\gt, Volume \thevolumenumber\ (\thevolumeyear)\hfill\fi}\vss
}}
\else
\makeatletter
\def\@oddhead{{\small\lhead\ifnum\count0=\startpage ISSN 1364-0380 (on line)
1465-3060 (printed) \hfill {\lnum\number\count0}\else\ifodd\count0
\def\\{ }\ifx\theshorttitle\relax \thetitle \else\theshorttitle\fi\hfill
{\lnum\number\count0}\else\def\\{ and }{\lnum\number\count0}
\hfill\ifx\theshortauthors\relax 
\theauthors\else\theshortauthors\fi\fi\fi}}\def\@evenhead{\@oddhead}
\def\@oddfoot{\small\lfoot\ifnum\count0=\startpage\copyright\ \gtp\hfill\else
\gt, Volume \thevolumenumber\ (\thevolumeyear)\hfill\fi}
\def\@evenfoot{\@oddfoot}
\makeatother
\fi

\newwrite\gtoutfile
\long\gdef\makeheadfile{  
{\def\\{, }\def\s{ }
\immediate\openout\gtoutfile head.xxx
\immediate\write\gtoutfile{To: math@arxiv.org}
\immediate\write\gtoutfile{Subject: put or rep NNNNN:pppp}
\immediate\write\gtoutfile{--text follows this line--}
\immediate\write\gtoutfile{Proxy-for: \ifx\theasciiauthors\relax
\theauthors\else\theasciiauthors\fi\s<\ifx\theasciiemail\relax\theemail\else\theasciiemail\fi>}
\immediate\write\gtoutfile{\noexpand\\}
\immediate\write\gtoutfile{Authors: \ifx\theasciiauthors\relax
\theauthors\else\theasciiauthors\fi}
{\def\\{ }\immediate\write\gtoutfile{Title: \ifx\theasciititle\relax
\thetitle\else\theasciititle\fi}}
\immediate\write\gtoutfile{Subj-class: GT or SG or MG etc}
\immediate\write\gtoutfile{MSC-class: \theprimaryclass\ifx\thesecondaryclass\relax\else, \thesecondaryclass\fi}
\immediate\write\gtoutfile{Journal-ref: Geom. Topol. \thevolumenumber
(\thevolumeyear) \startpage-\finishpage}
\immediate\write\gtoutfile{Comments: Published by Geometry and Topology at}
\immediate\write\gtoutfile{\s\s http://www.maths.warwick.ac.uk/gt/GTVol\thevolumenumber/paper\thepapernumber.abs.html}
\immediate\write\gtoutfile{\noexpand\\}
\immediate\write\gtoutfile{}
\ifx\theasciiabstract\relax
\immediate\write\gtoutfile{\theabstract}\else
\immediate\write\gtoutfile{\theasciiabstract}\fi
\immediate\write\gtoutfile{}
\immediate\write\gtoutfile{\noexpand\\}
\immediate\write\gtoutfile{}
\immediate\closeout\gtoutfile}}  

\def\maketitlepage{\maketitlep\makeheadfile}
\let\maketitle\maketitlepage

\def\ifplaintex{\expandafter\ifx\csname documentclass\endcsname\relax}

\ifplaintex 
\else
\fi

\expandafter\ifx\csname beginpicture\endcsname\relax
\expandafter\ifx\csname documentclass\endcsname\relax
\input pictex \else
\input prepictex \input pictex \input postpictex \fi\fi

\def\gt{{\mathsurround=0pt\it $\cal G\mskip-2mu$eometry \&\ 
$\cal T\!\!$opology}}        

\def\gtp{{\mathsurround=0pt\it $\cal G\mskip-2mu$eometry \&\ 
$\cal T\!\!$opology $\cal P\!$ublications}}  

\def\lognumber#1{\def\thelognumber{#1}}
\def\volumenumber#1{\def\thevolumenumber{#1}}
\def\papernumber#1{\def\thepapernumber{#1}}
\def\volumeyear#1{\def\thevolumeyear{#1}}

\def\pagenumbers#1#2{\def\startpage{#1}\def\finishpage{#2}}
\def\published#1{\def\publishdate{#1}}
\def\proposed#1{\def\theproposer{#1}}
\def\seconded#1{\def\theseconders{#1}}
\def\received#1{\def\receiveddate{#1}}
\def\revised#1{\def\reviseddate{#1}}
\def\accepted#1{\def\accepteddate{#1}}

\long\def\asciiabstract#1{\long\def\theasciiabstract{#1}}

\def\shorttitle#1{\def\theshorttitle{#1}}

\let\\\par\let\thelognumber\relax
\let\thevolumenumber\relax\let\thepapernumber\relax
\let\thevolumeyear\relax\let\thesamplenumber\relax\let\startpage\relax
\let\finishpage\relax\let\publishdate\relax\let\receiveddate\relax
\let\reviseddate\relax\let\accepteddate\relax\let\theasciititle\relax
\let\theasciiauthors\relax
\let\theasciiabstract\relax
\let\theasciiemail\relax\let\theshortauthors\relax\let\theshorttitle\relax

\long\def\maketitlep{   

\count0=\startpage

\gt\hfill      
\beginpicture
\setcoordinatesystem units <0.33truein, 0.33truein> point at 2.2 0.9
\setplotsymbol ({$\cal G$})
\plotsymbolspacing=9truept
\circulararc 315 degrees from 0 1 center at 0 0
\setplotsymbol ({$\cal T$})
\circulararc 315 degrees from 1 -1 center at 1 0
\endpicture
\break
{\small\ifx\thesamplenumber\relax 
Volume \else Sample
\fi\thevolumenumber\ (\thevolumeyear)
\startpage--\finishpage\nl
Published: \publishdate}
\vglue 0.5truein plus 0.4fil minus 0.1truein

{\parskip=0pt\leftskip 0pt plus 1fil\def\\{\par\smallskip}{\ifplaintex\large
\else\Large\fi\bf\thetitle}\par\medskip}   

\vglue 0pt plus 0.1fil 

{\parskip=0pt\leftskip 0pt plus 1fil\def\\{\par}{\sc\theauthors}
\par\medskip}

\vglue 0pt plus 0.1fil 

{\small\parskip=0pt\let\newline\\
{\leftskip 0pt plus 1fil\def\\{\par}{\sl\theaddress}\par}
\expandafter\ifx\theemail\relax    
\relax\else\vglue 5pt plus 0.02fil minus 2pt\def\\{\stdspace{\rm 
and}\stdspace} 
\cl{Email:\stdspace\tt\theemail}\fi
\ifx\theurl\relax                  
\relax\else\vglue 5pt plus 0.02fil minus 2pt\def\\{\stdspace{\rm 
and}\stdspace}
\cl{URL:\stdspace\tt\theurl}\fi\par}

\vglue 7pt plus 0.3fil minus 3pt

{\bf Abstract}
\vglue 5pt plus 0.1fil minus 2pt

\theabstract

\vglue 7pt plus 0.3fil minus 3pt

{\bf AMS Classification numbers}\quad Primary:\quad \theprimaryclass

Secondary:\quad \thesecondaryclass

\vglue 5pt plus 0.3fil minus 2pt

{\bf Keywords}\quad \thekeywords

\vglue 10pt plus 0.5fil minus 5pt

{\small  Proposed: \theproposer\hfill Received: \receiveddate\nl
Seconded: \theseconders\hfill 
\ifx\reviseddate\relax                         
Accepted: \accepteddate                        
\else
Revised: \reviseddate                          
\fi}
\eject
}       

\let\maketitlepage\maketitlep
\let\maketitle\maketitlepage

\font\phead=cmsl9 scaled 950
\font=cmsl9 scaled 1050
\font\pnum=cmbx10 scaled 913
\font\lnum=cmbx10 
\font\pfoot=cmsl9 scaled 950
\font=cmsl9 scaled 1050
\ifplaintex
\headline{\vbox to 0pt{\vskip -4.5mm\line{\small\phead\ifnum
\count0=\startpage ISSN 1364-0380 (on line)
1465-3060 (printed) \hfill {\pnum\folio}\else\ifodd\count0\def\\{ }%
\ifx\theshorttitle\relax\thetitle\else\theshorttitle\fi\hfill{\pnum\folio}
\else\def\\{ and }{\pnum\folio}\hfill\ifx\theshortauthors\relax\theauthors
\else\theshortauthors\fi\fi\fi}\vss}}
\footline{\vbox to 0pt{\vglue 0mm\line{\small\pfoot\ifnum\count0=\startpage
\copyright\ \gtp\hfill\else
\gt, Volume \thevolumenumber\ (\thevolumeyear)\hfill\fi}\vss
}}
\else
\makeatletter
\def\@oddhead{{\small\lhead\ifnum\count0=\startpage ISSN 1364-0380 (on line)
1465-3060 (printed) \hfill {\lnum\number\count0}\else\ifodd\count0
\def\\{ }\ifx\theshorttitle\relax \thetitle \else\theshorttitle\fi\hfill
{\lnum\number\count0}\else\def\\{ and }{\lnum\number\count0}
\hfill\ifx\theshortauthors\relax 
\theauthors\else\theshortauthors\fi\fi\fi}}\def\@evenhead{\@oddhead}
\def\@oddfoot{\small\lfoot\ifnum\count0=\startpage\copyright\ \gtp\hfill\else
\gt, Volume \thevolumenumber\ (\thevolumeyear)\hfill\fi}
\def\@evenfoot{\@oddfoot}
\makeatother
\fi

\newwrite\gtoutfile
\long\gdef\makeheadfile{  
{\def\\{, }\def\s{ }
\immediate\openout\gtoutfile head.xxx
\immediate\write\gtoutfile{To: math@arxiv.org}
\immediate\write\gtoutfile{Subject: put or rep NNNNN:pppp}
\immediate\write\gtoutfile{--text follows this line--}
\immediate\write\gtoutfile{Proxy-for: \ifx\theasciiauthors\relax
\theauthors\else\theasciiauthors\fi\s<\ifx\theasciiemail\relax\theemail\else\theasciiemail\fi>}
\immediate\write\gtoutfile{\noexpand\\}
\immediate\write\gtoutfile{Authors: \ifx\theasciiauthors\relax
\theauthors\else\theasciiauthors\fi}
{\def\\{ }\immediate\write\gtoutfile{Title: \ifx\theasciititle\relax
\thetitle\else\theasciititle\fi}}
\immediate\write\gtoutfile{Subj-class: GT or SG or MG etc}
\immediate\write\gtoutfile{MSC-class: \theprimaryclass\ifx\thesecondaryclass\relax\else, \thesecondaryclass\fi}
\immediate\write\gtoutfile{Journal-ref: Geom. Topol. \thevolumenumber
(\thevolumeyear) \startpage-\finishpage}
\immediate\write\gtoutfile{Comments: Published by Geometry and Topology at}
\immediate\write\gtoutfile{\s\s http://www.maths.warwick.ac.uk/gt/GTVol\thevolumenumber/paper\thepapernumber.abs.html}
\immediate\write\gtoutfile{\noexpand\\}
\immediate\write\gtoutfile{}
\ifx\theasciiabstract\relax
\immediate\write\gtoutfile{\theabstract}\else
\immediate\write\gtoutfile{\theasciiabstract}\fi
\immediate\write\gtoutfile{}
\immediate\write\gtoutfile{\noexpand\\}
\immediate\write\gtoutfile{}
\immediate\closeout\gtoutfile}}  

\def\maketitlepage{\maketitlep\makeheadfile}
\let\maketitle\maketitlepage

\lognumber{172}
\volumenumber{6}
\papernumber{12}
\volumeyear{2002}
\pagenumbers{361}{391}
\received{28 March 2001}
\revised{8 July 2002}
\accepted{9 July 2002}
\published{13 July 2002}
\proposed{Jean-Pierre Otal}
\seconded{Benson Farb, Walter Neumann}

\input rlepsf
\let\relabela\adjustrelabel

\newtheorem{fact}{Fact}
\newtheorem{sublemma}{Sub-lemma}

\newtheorem{claim}{Claim} 
\newtheorem{corollary}{Corollary}
\newtheorem{theorem}{Theorem}
\newtheorem{lemma}{Lemma}
\newtheorem{proposition}{Proposition}
\newtheorem{scholium}{Scholium}

\newtheorem{formula}{Formula}

\newcommand{\Bbb}[1]{%
{\bf #1}}
\newcommand{\dfn}[1]{%
{\it #1\/}}
\newcommand{\frak}[1]{%
{\bf #1}}

\newcommand{\tot}[1]{%
{ \left| \left| {#1}\right| \right|}}
\newcommand{\sech}{{\hbox{\rm sech}} }

\begin{document}

\title{Characterizing the Delaunay decompositions of\\compact 
hyperbolic surfaces}
\shorttitle{Characterizing Delaunay decompositions}

\author{Gregory Leibon}

\address{Hinman Box 6188, Dartmouth College\\Hanover  NH 03755, USA}

\email{leibon@dartmouth.edu}

\begin{abstract}
Given a Delaunay decomposition of a compact hyperbolic surface, one 
may record the topological data of the decomposition, together with 
the intersection angles between the ``empty disks'' circumscribing the
regions of the decomposition. The main result of this paper is a 
characterization of when a given topological decomposition and angle 
assignment can be realized as the data of an actual Delaunay
decomposition of a hyperbolic surface.\end{abstract}

\asciiabstract{
Given a Delaunay decomposition of a compact hyperbolic surface, one 
may record the topological data of the decomposition, together with 
the intersection angles between the `empty disks' circumscribing the
regions of the decomposition. The main result of this paper is a 
characterization of when a given topological decomposition and angle 
assignment can be realized as the data of an actual Delaunay
decomposition of a hyperbolic surface. 
}

\keywords{Delaunay triangulation, hyperbolic polyhedra, disk pattern}

\primaryclass{52C26}
\secondaryclass{30F10}
\maketitlepage

\begin{section}{Introduction}\label{sec:1}

Given a Delaunay decomposition of a compact hyperbolic surface, one
may record the topological data of the decomposition, together with
the intersection angles between the ``empty disks'' circumscribing the
regions of the decomposition.  The main result of this paper (Theorem
\ref{cir1}) is a characterization of when a given topological
decomposition and angle assignment can be realized as the data of an
actual Delaunay decomposition of a closed hyperbolic surface.  As a
consequence of this characterization, we get a characterization of the
convex ideal hyperbolic polyhedra associated to a compact surface with
genus greater than one (Corollary \ref{thurs}); this result is a
generalization of the convex ideal case of the Thurston--Andreev
Theorem. Corollary \ref{thurs} emerges naturally from Theorem
\ref{cir1}, because the main ingredient in exploring the Delaunay
decomposition is a triangulation production result (Proposition
\ref{tri1}) whose proof relies on properties of the volume of
hyperbolic polyhedra.

This paper is organized as follows.  In Section \ref{sec:2}, we
 formulate Theorem \ref{cir1} and reduce its proof to three
 fundamental propositions: Propositions \ref{tri1}, \ref{fund}, and
 \ref{emf}.  Section \ref{Hyp} contains the proof of Proposition
 \ref{fund} along with an exploration of Theorem \ref{cir1}'s
 relationship to hyperbolic polyhedra.  In Section \ref{sec:3}, we
 exploit properties of hyperbolic volume in order to prove Proposition
 \ref{tri1}.  Section \ref{sec:4} contains the proof of Proposition
 \ref{emf}, a linear programming problem.  Section \ref{sec:5}
 contains a discussion of consequences, generalizations, and natural
 questions.  In particular, Section \ref{sec:5} contains a discussion
 of how we may obtain from Theorem \ref{cir1} a natural polyhedral
 tessellation of the Teichm\"uller space of a compact Riemann surface
 with genus greater than one and at least one distinguished point.

\end{section}

\begin{section}{The Delaunay decomposition}\label{sec:2}

To begin this section we will recall how to construct the
Delaunay decomposition of a compact boundaryless hyperbolic surface, $G$,
with respect to  a finite set of
distinct specified points, $V=\{v_i\}_{i=1}^{|V|}$.

\medskip
{\bf Delaunay construction}\qua
As a first step, lift $G$ to its universal cover, 
the hyperbolic plane, which we will always denote as $H^2$. To the
inverse image of $V$, $\pi^{-1}\left( V \right)$,  we apply Delaunay's
empty sphere method, originally  introduced by Delaunay in \cite{De}.
Namely, if a  triple in $\pi^{-1}\left (V \right)$ lies on the
boundary of a  disk with interior   empty of
other points in $\pi^{-1}\left (V \right)$, then  look at all the points in
$\pi^{-1}\left (V \right)$ on this disk's  boundary and take the 
convex hull of
these points.  This procedure will tile $H^2$ with geodesic polygons 
and this tiling will be invariant under $G$'s deck group.
Let us call the resulting   ``polygonal decomposition"
 (to be defined  carefully in Section \ref{sec:5})  the {\dfn{Delaunay
decomposition}}  of $G$ with respect to $V$.  

In this section, we will set up some preliminary
 topological structures that will be used keep 
track of the combinatorics of such  decompositions. 
 For the sake of simplicity 
 we will assume that the  topological structure 
underlying our decomposition is
a triangulation and  restrict our attention to compact boundaryless
surfaces. Also, we assume that our geometry 
is always hyperbolic and, in particular, a geodesic triangulation will
always mean a geodesic triangulation of a hyperbolic  
surface.  In Section \ref{sec:5},  we will carefully articulate 
a more general, but technically less pleasant, context into which 
 the arguments explored in this paper will still apply.  This will 
include dealing  with the  above needed ``polygonal decompositions", as
well as allowing our hyperbolic surfaces to have 
boundary, singularities, and corners.

The following notation will be used throughout.  $\frak{T}$ will 
denote  a triangulation  and   $\tot{\frak{T}}$ will  denote 
$\frak{T}$'s total space.
$\frak{T}$'s cells will be referred to as follows: 
the set of vertices will be denoted as  $V$, the set of the edges will be 
denoted by $E$, and the set
of faces will be denoted by $F$. If $S$ is a collection of cells in 
 $\frak{T}$, then let
$\tot{S}$  denote the subspace of $\tot{\frak{T}}$ corresponding to it. We
use  this notation since $|K|$ will always mean the cardinality of a set
$K$. Also, we will abuse notation a bit and
let, when appropriate, $c \in S$ mean that $\tot{c} \subset \tot{S}$.

In order to state the main theorem, we will need to carefully
articulate the geometric quantities that arise in a Delaunay decomposition 
and, in particular, we will need to describe how to keep track of
the intersection angles between neighboring ``empty disks".

\begin{subsection}{Characterizing the Delaunay
decomposition}\label{DelDecSec}

Given a geodesic triangulation  and $t \in F$ we may embed $t$ 
in $H^2$.  The lifted  vertices of $t$ lie on the boundary of an 
Euclidean disk in the Poincare  disk model of $H^2$. 
The intersection of this disk and the Poincare
disk  will be called the polygon's {\dfn{circumscribing disk}}.
Notice, the boundary of the circumscribing disk can be intrinsically described by
either a hyperbolic circle, a horocircle, or a banana 
(a connected component of the set of points of a fixed distance from a
geodesic).  

We will now  carefully measure the angle between the
circumscribing disks of  neighboring faces.  
Fix a $t \in F$ that arises from a geodesic triangulation and embed $t$ 
in $H^2$.  For an  $e \in t$  let $\gamma^e$ be the   geodesic
in $H^2$ containing the edge corresponding to  $e$.
Let $\theta^{e}_{t}$ denote the angle
between the boundary of $t$'s
circumscribing disk and $\gamma^e$ as measured on
the side of $\gamma^e$ containing $t$'s embedded image.  
For any  $e \in E$ there are faces $t_1$ and
$t_2$ such that $e\in t_1$ and $e \in t_2$ and we will let the
{\dfn{intersection angle}}
at $e$ between
the circumscribing disks
containing $\tot{t_1}$ and $\tot{t_2}$ be defined as
$\theta^{e}_{t_1}+\theta^{e}_{t_2}$.

In order to combinatorially keep track of these angle values,
we first fix a triangulation, let $\{\psi_e\}_{e \in E}$ be a  
basis of $\Bbb R^{|E|}$, and  let  $\psi^e$ denote the dual
of  $\psi_e$ in  $\{\psi_e\}_{e \in E}$'s  dual basis. 
Throughout this paper,  $p \in \Bbb R^{|E|}$ 
will denote a vector in such a basis and
$\psi^e(p)$ will be called the {\dfn{informal
intersection angle complement}} at $e$. 
Furthermore, let
$\theta^e(p) = \pi - \psi^e(p)$ be called the {\dfn{informal
intersection angle}} at $e$. 
Notice, associated to  any geodesic triangulation we 
have a  $p\in \Bbb{R}^{|E|}$
determined by the actual intersection angles
between the circumscribing disks.   In this situation $p$
is said to be {\dfn{realized}} by this geodesic triangulation, which
will be called $p$'s {\dfn{realization}}. We will prove the following fact
 in Section \ref{sec:3}.

\begin{fact}\label{unique}
If $p\in \Bbb{R}^{|E|}$ is realized by a geodesic triangulation, 
then this realization is unique.
\end{fact}

A realization of $p \in \Bbb{R}^{|E|}$
forces a pair of  necessary linear constraints
on $p$.
The most important of these constraints
captures the fact that in an actual geodesic triangulation,  each 
face has positive area, see formula \ref{sett}. 
 Call  $p \in \Bbb R^{|E|}$ {\dfn{feasible}} if  for any set of
faces $S$ we have, that 
\[   \sum_{e \in S}  \theta^e(p) > \pi |S|. \]
In the context of a Delaunay decomposition, a
feasibility inequality  tends to an equality as 
vertices approach each other.
The other  necessary linear constraint on $p \in \Bbb R^{|E|}$ , in order
 for it to be realized, is that at each vertex
 the hyperbolic structure is non-singular.
Namely, 
if at  every vertex $v$, $p$ satisfies 
\[ \sum_{\{e \mid  v \in e\}} \psi^e(p) = 2 \pi, \]
 then  we will call  $p$ {\dfn{non-singular}}.

In order to prove the existence of a realization for a given  
$p\in \Bbb{R}^{|E|}$,
we will assume that $p$ 
 arises  not only from  a geodesic triangulation but, in fact, from a 
Delaunay decomposition.  
In an actual Delaunay decomposition, it is straight forward 
to check that at each edge $\psi^e(p) \in (0,\pi)$, hence, we will call 
$p$  {\dfn{Delaunay}} if $\psi^e(p) \in (0,\pi)$ for all $e$.
The following is our main theorem.  

\begin{theorem}\label{cir1}
A non-singular Delaunay $p
\in \Bbb R^{|E|}$ is realized by a unique geo\-desic triangulation 
if and only if $p$ is feasible.  
\end{theorem}

{\bf Context}\qua
A Euclidean version of this theorem   can  be found in Bowditch's
\cite{Bo},
where it is proved using techniques similar to those found in Thurston's
proof
of the Thurston--Andreev Theorem in \cite{Th}. As with Bowditch's result,
 the non-singularity condition is unnecessary and,  in Section \ref{sec:5}, we will drop it.    The proof here relies on
the triangulation production result, Proposition \ref{tri1}, stated in
Section
\ref{san}.

\end{subsection}
 
\begin{subsection}{Triangulation Production}\label{san}

In a
triangulation there are exactly  $3|F|$ slots in which one can
insert possible triangle angles, and we will identify the possible triangle 
angle values  with the coordinates  of 
$\Bbb R^{3|F|}$.  These coordinates can be naturally indexed by 
$(e,t) \in E \times F$ by identifying the angle slot in $t$ opposite to 
$e$ with a basis vector $\alpha^t_e$.
Throughout this paper,  $x \in \Bbb R^{3|F|}$ will denote a vector 
in this basis and $\alpha^e_t$ will denote the dual of $\alpha^t_e$ in 
the dual basis.
In particular,  $\alpha^e_t(x)$ will represent the angle 
in the $(e,t)$ angle slot.
Usually we will view these vectors and covectors geometrically.  A vector will be
expressed by placing its coefficients in a copy of the triangulation  with
dashed lines and a covector will contain its   coefficients  
in a copy of the triangulation with solid
lines, see figure \ref{cv}. 
Angle slots not pictured
will always
be assumed to have  zero as their coordinate's value.  Notice, the pairing of a
vector and a covector is achieved  by placing the copy of the  
 triangulation corresponding to the covector on top of the  triangulation 
corresponding to the vector and multiplying  the numbers living in the
same angle slots.  In figure \ref{cv},  we see the covectors and vectors
that we  will be needing. For easy
reference, we will now sum up the relevant definitions related 
directly to the
pictured covectors and vectors.

\medskip 
{\bf Definition 1}\qua
Assume $x \in \Bbb{R}^{3|F|}$.  For a
triangle $t$ 
with edges $\{e_1,e_2,e_3\}$, 
let $d^t (x)   =
\{\alpha^{e_1}_t(x),\alpha^{e_2}_t(x),\alpha^{e_3}_t(x) \}$
and call
$d^t(x)$  the {\dfn{angle data}} associated to $t$.
$k^t(x) =
\sigma^t(x)-\pi$ will be called $t$'s {\dfn{curvature}}. 
Notice, if $k^t(x) < 0$ , then we may form an actual hyperbolic triangle
with the
angles in
$d^t(x)$ which will be called $d^t(x)$'s
{\dfn{realization}}.   For each
$e
\in t$ denote the length of the edge
$e$ with respect  to $d^t(x)$'s
realization as $l^e_t(x)$.
For an  edge
$e$  
sharing the  triangles
$t_1$ and $t_2$ we will let
\[  \psi^{e} =
\psi^{e}_{t_1} + \psi^{e}_{t_2},\]
with $\psi^e_t$ the covector defined in
figure one.
We will call  $\psi^e$   the {\dfn{informal angle
complement}}
at $e$. The {\dfn{informal intersection angle}}  is  defined as
\[
\theta^e(x) = \pi - \psi^e(x).\] 
Call $x$ non-singular.
Let an {\dfn{angle system}}  be a point in
  \[
\frak{N} = \{ x  \mid k^t(x) < 0 \ \forall t, 
\alpha^e_t(x) \in (0,
\pi) \ \forall (e,t)
, r^v(x)=2\pi  \,\forall v \}. \]
A {\dfn{conformal
deformation}}
 will be a vector in
\[ \frak{C} = span\{ w_e \mid \,\forall e
\},\]
and we will call $x$ and $y$  {\dfn{conformally equivalent}} if $x -
y \in
\frak{C}$.
Let
\[\frak{N}_x = (x + \frak{C}) \bigcap \frak{N}\]
be
called the {\dfn{conformal class of}} $x$.
 Let
\[ \Psi\co  \Bbb R^{3|F|}
\rightarrow \Bbb R^{|E|} \]
be the linear mapping  given by
\[ \Psi(x) =
\sum_{e \in E} \psi^{e}(x) \psi_e.\]
We call $x$ {\dfn{Delaunay}}   or {\dfn{non-singular}} if  $\Psi(x)$  is.

\begin{figure}[ht!]
\cl{\relabelbox\small\let\ss\scriptsize
\let\sss\scriptscriptstyle
\epsfxsize 3in \epsfbox{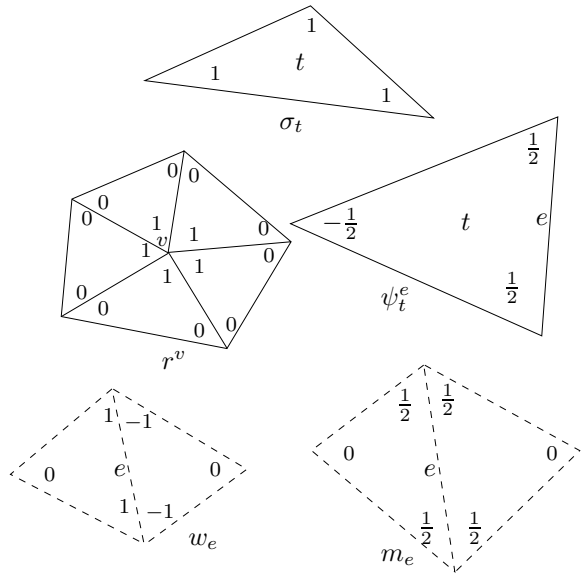}
\relabel {12}{$1\over2$}
\relabel {121}{$1\over2$}
\relabel {122}{$1\over2$}
\relabel {123}{$1\over2$}
\relabel {124}{$1\over2$}
\relabel {125}{$1\over2$}
\relabel {-12}{$-{1\over2}$}
\relabel {0}{\ss$0$}
\relabel {01}{\ss$0$}
\relabel {02}{\ss$0$}
\relabel {03}{\ss$0$}
\relabel {04}{\ss$0$}
\relabel {05}{\ss$0$}
\relabel {06}{\ss$0$}
\relabel {07}{\ss$0$}
\relabel {08}{\ss$0$}
\relabel {09}{\ss$0$}
\relabel {010}{\ss$0$}
\relabel {011}{\ss$0$}
\relabel {012}{\ss$0$}
\relabela <-4pt,0pt> {013}{\ss$0$}
\relabel {1}{\ss$1$}
\relabel {1a}{\ss$1$}
\relabel {1b}{\ss$1$}
\relabel {1c}{\ss$1$}
\relabel {1d}{\ss$1$}
\relabel {1e}{\ss$1$}
\relabel {1f}{\ss$1$}
\relabel {1g}{\ss$1$}
\relabel {1h}{\ss$1$}
\relabela <-2pt,0pt> {1i}{\ss$1$}
\relabel {m1}{\ss$-1$}
\relabel {mm1}{\ss$-1$}
\relabel {v}{\ss$v$}
\relabel {y}{$\psi^e_t$}
\relabel {s}{$\sigma_t$}
\relabel {r}{$r^v$}
\relabel {w}{$w_e$}
\relabel {m}{$m_e$}
\relabel {e}{$e$}
\relabel {e1}{$e$}
\relabel {e2}{$e$}
\relabela <0pt,-2pt> {tt}{$t$}
\relabela <0pt,-2pt> {ttt}{$t$}
\endrelabelbox}
\caption{Here are the relevant vectors
and covectors.
}
\label{cv}
\end{figure}

By
the Gauss--Bonnet Theorem, we know that the curvature, $k^t(x)$,  is  the
curvature in  a geodesic triangle with  angle
data $d^t(x)$.  Angle
systems
are precisely the $x \in \Bbb{R}^{3|F|}$  where the curvature
is
negative
and all angles are realistic. In particular, the actual angle
data of
 a geodesic triangulation is a special case of an angle
system. 
If  $u \in \Bbb{R}^{3|F|}$ 
is the
data of a geodesic triangulation, then we will call  an $u$
{\dfn{uniform}}.

Our goal is to take a point in
$\frak{N}$  and  deform it into a uniform
point.   Such deformations are
located in the  affine space of 
conformal deformations,
$\frak{C}$.  The
following is a triangulation production result
whose proof (accomplished
in Section \ref{sec:3})  
will  utilize this conformal deformation
strategy.

\begin{proposition}\label{tri1}
Every non-singular Delaunay
angle system is conformally\break equivalent
to a unique uniform angle
system.
\end{proposition}

{\bf Context}\qua
The proof presented here  of
this triangulation  production result  relies on
the use an objective
function.   It should be acknowledged that, the  use of an
objective
function in the setting of triangulation production  was first
accomplished by Colin de Verdi\'ere in \cite{Co}. Our objective function
turns out to be   related in a rather magical way to hyperbolic volume,
and 
that a connection  between triangulation production and hyperbolic
volume 
should exist has its origins  in  Br\"agger's beautiful paper
\cite{Be}. 
The particular volume exploited here was  observed by my PhD
thesis advisor,  Peter
Doyle. In the Euclidean case, this technique was
carried out by  Rivin in 
\cite{Ri2}.

The relationship of Theorem
\ref{cir1}  to Proposition \ref{tri1}  comes from
the following
fundamental observation, which will be proved
in Section
\ref{Hyp}.

\begin{proposition}\label{fund}
Let $u$ be a  uniform angle
system and let  $e \in E$,  then
$\theta^e(u)$ is  the intersection angle
of the circumscribing circles of the
hyperbolic triangles sharing
$e$.
\end{proposition}

Proposition \ref{fund} immediately shows us that
the use of the  notation  
$\psi^e$ and $\theta^e$, in section \ref{DelDecSec},
is compatible with the use in 
definition  1, namely, for a uniform angle
system $u$ we have that 
$\theta^e(u)$ really is the intersection angle
$\theta^e$.

Much of what takes place here relies
on certain basic
invariants of conformal  deformations.  For a simple example
of a conformal invariant,  notice  that,  $r^v(x)$ (from figure \ref{cv})
satisfies
\[ r^v(y) = r^v\left(x + \sum_{e \in E} B^{e} w_{e}\right) 
= r^v(x) + \sum_{e \in E} B^{e} r^v(w_{e} )  =  r^v(x) ,\]
which tells us that non-singularity is a  conformal invariant. 
Similarly we have that 

\vspace{-5mm}
\[ \psi^e(y) = \psi^e\left(x + \sum_{g \in F} B^g w_{g}\right) 
= \psi^e(x) + \sum_{g \in F} B^g \psi^e(w_{g} ) 
= \psi^e(x), \]
and, hence,
$\psi^e$ and $\theta^e$ are also conformal invariants.
In fact, this preservation of the formal intersection angle
is why  these   transformations are called conformal (see
\cite{Le1} for a deeper reason). 
The next lemma expresses the fact that, the elements of 
$\frak{C}$ are the only conformal invariants in this sense.

\begin{lemma}\label{idpun}
For any $p \in \Bbb R^{|E|}$, $\Psi^{-1}(p) = \frak{C} + x$ for 
any $x \in \Psi^{-1}(p)$.
\end{lemma} 
\proof
Since the pairing of $\psi^{e}$ with 
the vector $m_{e}$,  in figure \ref{cv}, 
satisfies   $\Psi(m_{e})\! = \psi_e$ for each edge $e$, $\Psi$ has rank $E$ .  
In particular, the null space of $\Psi$ is  $3F - E  =E$
dimensional.
However, by the conformal
invariance of $ \psi^e $, the null space contains the 
$E$ dimensional space $\frak{C}$; hence, the null space of $\Psi$
is  $\frak{C}$.
\endproof

Let  $ \frak{N}^p= \frak{N}_x $
for some $x \in \Psi^{-1}(p)$, and note,
as a consequence of Lemma \ref{idpun}, Proposition \ref{tri1}, 
and Proposition \ref{fund}, that a non-singular Delaunay 
$p\in {\Bbb R}^{|E|}$ is realized by  a unique geodesic triangulation if and
only if $ \frak{N}^p$ is non-empty.
Hence, Theorem \ref{cir1} is now seen to be implied by the fact that
$r^v(x)$ is conformally invariant and the following proposition, 
which is proved  in Section \ref{sec:4}.

\begin{proposition}\label{emf}
For a Delaunay $p \in \Bbb R^{|E|}$, we have that 
$\frak{N}^p$ is non-empty if and only if 
$p$ is feasible.
\end{proposition}

{\bf Context}\qua
Proposition \ref{emf} can be phrased as
a linear programming question. 
In particular, since we know a solution exists, the simplex algorithm
can be used to actually find the solution, 
see \cite{Dan}.  Knowing that a solution exists is
usually referred to as a feasibility criteria and is the motivation for the 
terminology feasible.
Similar geometrically motivated problems arose in
\cite{Co}  and \cite{Ri3}, where they were solved using 
techniques from linear programming.
These techniques appear to the author be much 
more complex to implement in this setting, hence, the
 optimization strategy described in Section \ref{sec:4}.

 \end{subsection}

\end{section}

\begin{section}{Hyperbolic polyhedra: The proof of
 Proposition \ref{fund}}\label{Hyp}

Let us call $H^2$ the hyperbolic plane in $H^3$ viewed as  in figure \ref{hit}.  
The inversion, $I$,  through the sphere of radius $\sqrt{2}$
centered at the south pole interchanges our specified $H^2$ with the upper
half of  the sphere at infinity, $S^{\infty}_u$.  Notice, when viewed
geometrically,  this map sends a point $p \in H^2$ to the point where
 the geodesic
perpendicular to $H^2$ containing $p$  hits $S^{\infty}_u$ (see figure
\ref{hit}).  In particular,
 being an inversion,  any circle in the $xy$--plane is
mapped to a circle on the sphere at $\infty$, $S^{\infty}$.
The use of this mapping will require the construction of an object that
will be  crucial in proving Proposition \ref{tri1}.

\begin{figure}[ht!]
\cl{\relabelbox\small\epsfxsize 2in \epsfbox{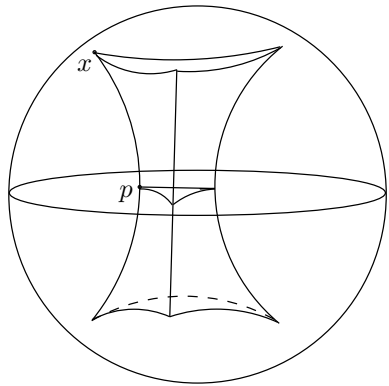}
\relabel {p}{$p$}
\relabel {x}{$x$}
\endrelabelbox}
\caption{\label{hit} We have our specified hyperbolic plane
$ H^2 \subset H^3$  realized as the intersection of the unit sphere at the
origin with the $xy$--plane in ${\bf R}^3$ via the Poincare  disk model of
$H^3$.  The point $p$, in the figure,  is mapped  to the point
labeled $x$ under the inversion $I$.   We also are viewing a triangle on
that plane and its associated prism in this model.}
\end{figure}

\medskip
{\bf Prism construction}\qua
Place a hyperbolic triangle on a copy of $H^2 \subset H^3$. 
Let its associated {\dfn{prism}} be the convex hull of 
the  set  consisting of the triangle 
unioned with  the geodesics 
perpendicular to this $H^2 \subset H^3$ 
going through the triangle's
vertices, as visualized in figure \ref{hit}.  
Denote the prism  relative to the hyperbolic triangle
constructed from the data $d = \{A,B,C\}$ as 
$P(d)$.

Now back to our proof.  Let $u$ be a uniform angle system, and let $t_1$
and $t_2$ be a pair of triangles sharing the edge $e$.  Place them next to
each other in our $H^2$ from figure \ref{hit}.  Notice,
$t_1$ and $t_2$ have circumscribing  circles in the $xy$--plane, which
correspond to either circles, horocircles, or bananas in $H^2$. 
Since the Poincare model is conformal, the intersection angle of these
circles is precisely the hyperbolic intersection angle.   
Being an inversion, $I$ is conformal.  In particular, these circles are sent
to  circles at infinity intersecting at the same angle  and going
through the ideal points of the neighboring $P(d_{t_1}(u))$ and
$P(d_{t_2}(u))$.   But these  circles at infinity  are also the intersection  of $S^{\infty}$ with the spheres representing 
the hyperbolic planes forming the top faces of 
$P(d_{t_1}(u))$ and $P(d_{t_2}(u))$.   
So the intersection angle of these spheres is precisely the sum
of the angles inside $P(d_{t_1}(u))$ and $P(d_{t_2}(u))$ at the edge
corresponding to $e$,   which we will now see is $\theta^e(u)$.  
In fact, we will show that this geometric decomposition of the
intersection angle is precisely the decomposition
 \[ \theta^e(u)= \theta^e_{t_1}(u) + \theta^e_{t_2}(u).\]

\begin{figure}[ht!]
\cl{\relabelbox\small
\epsfxsize 1.7in \epsfbox{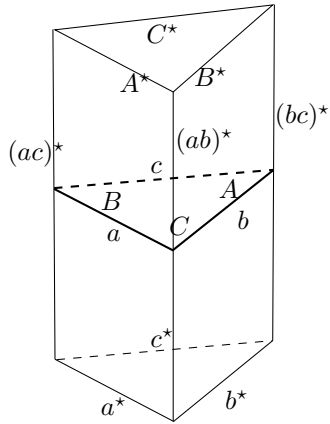}  
\relabel {A1}{$A$}
\relabel {B1}{$B$}
\relabel {C1}{$C$}
\relabel {A2}{$A^\star$}
\relabel {B2}{$B^\star$}
\relabel {C2}{$C^\star$}
\relabel {a1}{$a$}
\relabel {b1}{$b$}
\relabel {c1}{$c$}
\relabel {a2}{$a^\star$}
\relabel {b2}{$b^\star$}
\relabel {c2}{$c^\star$}
\relabel {ac}{$(ac)^\star$}
\relabel {ab}{$(ab)^\star$}
\relabel {bc}{$(bc)^\star$}
\endrelabelbox}
\caption{\label{clam} The notation for an ideal prism associated to 
hyperbolic angle data $\{A,B,C\}$, viewed for convenience  in the Klein
model.  }
\end{figure}

Now pick an $i$ and let  $d^{t_i}(u)$ be denoted  by $\{A,B,C\}$.  Assume
our specified edge
$e$ corresponds to the
$a$ in figure \ref{clam}.
From figure \ref{horo1}, we see the angles in figure \ref{clam} satisfy the
system of linear equation telling us that interior angles of the prism sum
to
$\pi$ at each vertex of the prism.  Solving this system   for the needed
angle, 
$A^{\star}$, we find that indeed
\[ A^{\star} =\frac{\pi + A -B -C}{2 } = \theta^e_{t_i}(u), \]
as claimed.

\begin{figure}[ht!]
\cl{\relabelbox\small
\epsfxsize 3in \epsfbox{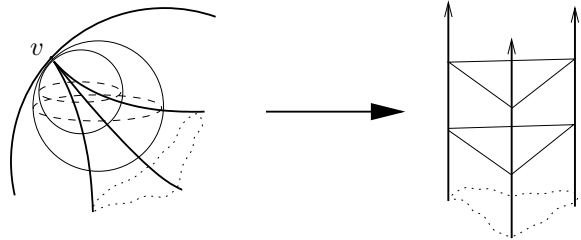}
\relabela <-2pt,2pt> {v}{$v$}
\endrelabelbox}  
\caption{\label{horo1} In this figure, we have  cut off an ideal
vertex of a convex hyperbolic polyhedron with a horoball. Horospheres are
flat planes in $H^3$ and our cut produces an Euclidean polygon in
this flat plane.  In particular, the interior angles at an ideal vertex of a
convex hyperbolic  polyhedron are those of an Euclidean polygon.  As
indicated, this is best  seen in the upper half space model of
$H^3$, with the ideal vertex being sent to
$\infty$.  (Nice applications of this observation can be found in \cite{Th}.) }
\end{figure}

We have now completed our proof of observation \ref{fund}. 
The ideas in this proof
 allow us to see how a geodesic triangulation may be 
used to  construct certain ideal hyperbolic polyhedra.
Given any geodesic triangulation  of $H^2 \subset H^3$, we may
form the  polyhedron 
$\bigcup_{t \in {\frak{T}}} P(d_t(u))$. We shall refer to the
polyhedra that can be formed via this construction applied to the lift
of a geodesic triangulation of a hyperbolic surface as
the  {\dfn{polyhedra associated
to surfaces of genus greater than 1}}. 
As a scholium to observation 
\ref{fund}, the dihedral angles in these polyhedra are precisely the
intersection angles  between the  circumscribing disks of
 neighboring faces in the triangulation.
In particular, the collection of 
dihedral angle complements can be naturally identified
with our  $p \in {\Bbb R}^{|E|}$, and results about
geodesic triangulations can  be translated into results
about such polyhedra. For example, Fact \ref{unique} tells us 
that such polyhedra are
uniquely determined by their dihedral angles, and we have the  
following corollary to Theorem \ref{cir1} concerning
the existence of such polyhedra.

\begin{corollary}\label{thurs}
A non-singular Delaunay $p
\in \Bbb R^{|E|}$ is realizable by a unique 
polyhedron associated
to surface of genus greater than 1 
if and only if $p$ is feasible.  
\end{corollary}

Notice,  the constraint that $\theta^e(x)$ is  in $(0,\pi)$, which
 corresponded to the Delaunay constraint,
corresponds to the  convex case among these polyhedra.
 As such, Corollary  \ref{thurs}  can be viewed as a characterization 
of the possible dihedral  angles in such convex
 polyhedra, hence, as a generalization of the  convex ideal  case of
the
 Thurston--Andreev Theorem.

\end{section}

\begin{section}{Hyperbolic volume: 
The proof of Proposition \ref{tri1}}\label{sec:3}

In this section we prove Proposition \ref{tri1} using an  argument dependent on a rather remarkable link between
hyperbolic volume and uniform structures.
Let the volume of the prism $P(d^t(x))$
be denoted $V(d^t(x))$.  We will be exploring the objective function  
\[H(y) = \sum_{t \in \frak{T}} V(d^t(y)) .\]
$H$'s behavior is at the heart of our arguments. 
Namely, Proposition \ref{tri1} follows immediately from 
the following claim about $H$'s behavior, while 
Fact \ref{unique} follows from Lemma 
\ref{idpun}  and the following claim.

\begin{claim}\label{tc}
View $H$ as a function on $N_x$ for some angle system $x$.
\begin{enumerate}   
\item
 $y$ is a critical point of $H$  if and only if $y$ is uniform.
\item
If $H$ has a critical point, then this critical point is unique.
\item
If $x$ is Delaunay,  then $H$ has a critical point.
\end{enumerate}
\end{claim}
\proof  
To prove the first part of the Claim \ref{tc} we compute   $H$'s differential. 
As usual, for  a function on a linear space like $\Bbb R^{3|F|}$,
we identify the tangent and cotangent spaces at
every point with  $\Bbb R^{3|F|}$  and  $(\Bbb R^{3|F|})^*$, and express our
differentials in the chosen basis.   In Section \ref{diff}, we prove the
following lemma.

\begin{lemma} \label{diff}
$dH(z) = \sum_{t \in F} \left( \sum_{e \in t} h^e_t(z)  \theta^e_{t}
\right)$ with the property that
$h^e_{t}(z)$ uni\-quely determines the length $l^e_{t}(z)$ (see Definition 1).
\end{lemma}
Recall that  $T_z(\frak{N}_x)$ may be identified with the vector
space of conformal deformations 
as described in definition 1.  Now, simply observe that $\theta^f_t(w_e) =
\pm \delta^f_e$, with the sign depending on whether $t$ contains the
negative or positive half of
$w_f$. So, at
a critical point $y$ we have  
\[ 0  =  \left| dH(y)(w_e) \right| = \left| h^e_{t_1}(y) -
h^e_{t_2}(y) \right| \]
where $t_1$ and $t_2$ are the faces sharing $e$.  From the  lemma \ref{diff}
and the fact that the $w_e$  span $T_y(\frak{N}_x)$,
 we have that
$l^e_{t_1}(y) = l^e_{t_2}(y)$ for each edge is equivalent
to  $y$ being critical.  Now, we finish off the proof of the first
part of Claim
\ref{tc}  by noting that, an angle system $u$  is uniform  if
and only if it fits together in the sense that the
 lengths  $l^e_{t_1}(u)$ and  
$l^e_{t_2}(u)$ are equal whenever $t_1$ and $t_2$ share an edge $e$.

In order to prove the second part of Claim \ref{tc},
 we need the following lemma, proved in
Section \ref{con}. 

\begin{lemma}\label{cony}
$H$   is a strictly concave,  smooth function  on $\frak{N}_x$ and 
continuous on
$\frak{N}_x$'s closure.
\end{lemma}
Now, a smooth concave function,like $H$,  has at most 
one critical point in any open
convex set, which  proves part 2.
To prove the third part of Claim \ref{tc}, it is useful to 
first isolate the technical property that makes 
Delaunay angle systems manageable.
We capture this property with the following definition.
 Call  $x \in \partial \frak{N}$ {\dfn{foldable}}  if $d^t(x)  =  \{
0,0,\pi \}$
for every $t$ where
either $k^t(x) =0$ or $d^t(x)  =
\{A,B,0\}$. Call a set in $\Bbb R^{|E|}$ {\dfn{unfoldable}} if the
set intersects $\frak{N}$ but contains no foldable points. 
We have the following lemma.

\begin{lemma} \label{delun}
A non-singular  Delaunay angle system is unfoldable.
\end{lemma}  
\proof
Let $x$ be a non-singular Delaunay angle system.
In this situation, unfoldability  is quite easy to verify since
$d^t(x + c)$ can never equal $\{0,0,\pi\}$ when $x +c \in \partial \frak{N}$
and  $\frak{c} \in \frak{C}$.
To see this  fact, assume to 
the contrary that for some $t$  and $\frak{c}$ we have 
$d^t(x + \frak{c}) = \{0,0,\pi\}$. 
Let $e$ be the  edge of $t$ across from
$t$'s  $\pi$
and let $t_1$ be $t$'s neighbor.
Now notice, if $x \in \frak{N}$ and
  $d^{t}(x) = \{A,B,C\} $ then,
since 
$B +C \leq A+B+C = \sigma^t(x) < \pi$ and $A<\pi$, we have  
\[ -\frac{\pi}{2} < -\frac{A}{2} <  
\psi^e_{t} = \frac{B+C-A}{2}
<  \frac{B +C}{2} < \frac{\pi}{2}.\]
In other words, for any  $x \in \frak{N}$  we have that
$\psi_{t}^e(x) \in \left(\frac{-\pi}{2},\frac{\pi}{2} \right)$.
In particular,  the conformally invariant $\psi^e(x) \in (0,\pi)$
 would have  to satisfy the inequality 
$\psi^e(x+c) = - \frac{\pi}{2} + \psi_{t_1}^e \leq 0 $,
contradicting the fact that $x$ is Delaunay.
\endproof

Now we will 
finish off the third part of Claim \ref{tc} by proving 
that, if $x+\frak{C}$ is unfoldable, then $H$ has a critical point.
To prove this, first notice, $\frak{N}_x$ is a
pre-compact open set and $H$ is a differentiable function
on $\frak{N}_x$ which is continuous in  $\frak{N}_x$'s closure, so 
we are assured of a
critical point if $H$'s maximum is not on $\frak{N}_x's$'s boundary. 
Notice, since
$\frak{N}_x$ is convex and $x+\frak{C}$ is unfoldable,
  that  for every  boundary
point, $y_0$, 
there is a direction $v$  such that when letting $l(s) = y_0 + s v$ 
we have that
$l\left([0,\infty)\right)$ is unfoldable.  If, for some $c >0$, we knew
that 
\[ \lim_{s \rightarrow 0^+} \frac{d}{ds}H(l(s)) \mbox{  } > \mbox{  } c, \]
then there would be some $\epsilon >0$ such that $H(l(s))$ is continuous and 
increasing on $[0,\epsilon)$ with $l((0,\epsilon)) \subset \frak{N}_x$ and,
in particular, $y_0$ could not have  been a point where$H$ achieved its maximum.  This fact is  immediately implied by  the following
lemma to be proved in Section \ref{bound}.

\begin{lemma}\label{pushin}
For every  point $y_0$ in $\partial \frak{N}$   and every direction $v$ 
such that  $l\left([0,\infty)\right)$ is unfoldable,   we have  
\[ \lim_{s \rightarrow 0^+} \frac{d}{ds}H(l(s)) \mbox{  } = \mbox{  } \infty. \]
\end{lemma} 
So, we indeed find that $H$ achieves its unique critical in $\frak{N}_x$, as
needed to complete the proof of   the claim.
\endproof

\begin{subsection}{The differential:  The computation of Lemma
\ref{diff}}\label{difff}

In this section we gain our needed understanding of the
differential as expressed in Lemma \ref{diff}.  
To get started, note the sum in $dH$ is over all faces, but
the fact concerns only each individual one.  So we may
restrict our attention to one triangle.
One way to prove Lemma \ref{diff}  
is to explicitly compute a formula for the volume
in terms of the Lobachevsky function and then find its
differential.  This method can be found carried out in
\cite{Le}.  Here we present an argument using Schlafli's formula for
 volume deformation.  This technique  has a wider
range of application as well as  being considerably more
interesting.

To start with, we will recall Schlafli's formula for a
differentiable family of compact convex polyhedra  with
fixed combinatorics. Let ${\bf{E}}$ denote the set of
edges and
$l(e)$ and $\theta(e)$  be the length and dihedral angle
functions associated to an edge $e$. Schlafli's formula is 
the following formula for the deformation  of the volume
within  this family

\[ dV =- \frac{1}{2} \sum_{e \in {\bf{E}}} l(e) d(\theta(e)). \]

In the finite volume case, when there  are ideal
vertices, the formula changes from measuring the length
of edges $l(e)$ to measuring the length of the 
cut off
edges $l^{\star}(e)$, a fact observed by Milnor  (see
\cite{Ri2} for a proof).  Let us now recall  how
$l^{\star} (e)$ is computed.  First 
fix a horosphere at each ideal vertex. Then note from any
horosphere to a point and between any pair of
horospheres there is a unique (potentially degenerate)
geodesic segment perpendicular  to the horosphere(s).  
$l^{\star}(e)$ is the signed length of this geodesic
segment;  it is given a positive sign if the geodesic is outside
the horosphere(s) and a negative sign if not. 
Schlafli's formula is
independent of the horosphere  choices in this
construction, and  I will refer to this fact as the
horoball independence  principle. It is   worth recalling the
reasoning behind this principle, since the ideas involved 
will come into play at several points in what follows.  

\medskip
{\bf The horosphere independence reasoning}\qua  Recall from the proof of observation 
\ref{fund}  that, at an ideal vertex $v$,  we have the sum of the dihedral  angles satisfying 
$ \sum_{e \in v} \theta(e) =  (n-2) \pi$, where $\{e \in v\}$ is the set
of edges containing $v$.
In particular 
\[ \sum_{e \in v} d \theta(e) =  0. \]
Looking at figure \ref{horo1}, we see by 
changing the horosphere at the ideal vertex $v$  that 
$l^{\star}(e)$  becomes  $l^{\star}(e) + c$ for
each $e \in v$  with $c$ a fixed constant.  Hence, by our
observation  about the angle differentials
\[-2dV =   \sum_{e \in {\bf E}} l^{\star}(e) d
\theta(e)= \sum_{e \in {\{e \in v\}^c} }l^{\star}(e) d \theta(e) + 
\sum_{e
\in v} (l^{\star}(e) + c )d \theta(e)  \]
and $dV$ is seen to be independent of the 
horosphere choices.
\endproof

Now let us look at our prism. Let the
notation for the cut off edge lengths coincide with the
edge names in figure
\ref{clam}.
 Since we may choose any horospheres, let us choose
those tangent to the hyperbolic plane which our prism is
symmetric across. In this case, note the lengths of 
$(ab)^{\star}$, $(bc)^{\star}$  and $(ac)^{\star}$ are zero.
Recalling from the proof of observation \ref{fund} that  
\[ A^{\star} =  \frac{\pi + A -B -C}{2} \]
and viewing $V(d^t(x))$ as a function on 
\[ \{(A,B,C) \in (0,\pi)^3 : 0<A+B+C<\pi \} \] 
we see from Schlafli's formula that
\[ dV =  - a^{\star} dA^{\star}  - b^{\star} d B^{\star}
 - c^{\star} d C^{\star} .\]

\begin{figure}[ht!]
\cl{\relabelbox\small
\epsfxsize 2in \epsfbox{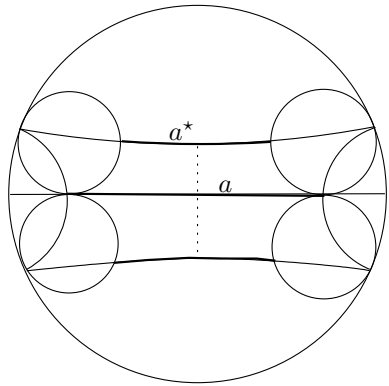}
\relabel {a}{$a$}
\relabel {a*}{$a^\star$}
\endrelabelbox}  
\caption{\label{facesl} Here we see a face of our prism
 containing $a$.  Also pictured are the horocircle slices of the
horospheres tangent to the hyperbolic plane through which our prism
is symmetric. }
\end{figure}

Note that Lemma \ref{diff} will follow  from the
following formula.

\begin{formula} 
\[a^{\star} = 2
\ln\left(\sinh\left(\frac{a}{2}\right)\right).
\]
\end{formula}\proof To begin this computation, look at the face of the
prism containing$a$
 as in figure \ref{facesl}.
 Notice, this face is
decomposed into four congruent  
quadrilaterals, one of which as been  as in figure
\ref{quad}. 
Note that, just as with the above reasoning
concerning the independence of horosphere choice, we have
an independence of horocircle choice and
\[\frac{a^{\star}}{2} = (t^{\star} - h^{\star}) - (h^{\star}
- s^{\star}) .\]
In fact, $t^{\star} - h^{\star}$  and  $h^{\star}
- s^{\star}$ are independent of this horocircle
choice as well, and it is these quantities we shall
compute.  
 
\begin{figure}[ht!]
\cl{\relabelbox\small
\epsfxsize 2.2in \epsfbox{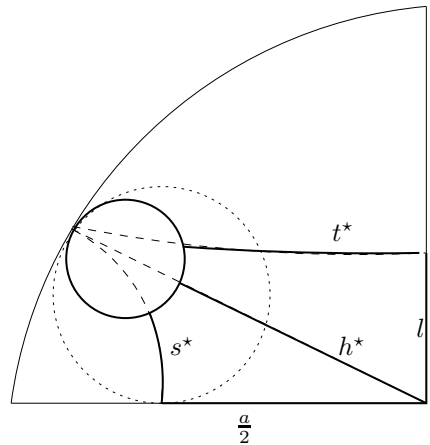}
\relabel {t}{$t^\star$}
\relabel {s}{$s^\star$}
\relabel {h}{$h^\star$}
\relabel {l}{$l$}
\relabel {a}{$a\over2$}
\endrelabelbox}  
\caption{Pictured here is one of the four triply right-angled
 quadrilaterals from figure \ref{facesl}.  Such quadrilaterals are
 known as Lambert quadrilaterals, and it is a well know relationship
 that $\tanh\left(\frac{a}{2}\right) = \sech(l)$, see for example
 \cite{Gr}.  (In fact it follows nearly immediately from the perhaps
 better known Bolyia--Lobachevsky formula.)}
\label{quad}
\end{figure}

Look at the figure \ref{quad} and notice, using the
 horocircle tangent to the $\frac{a}{2}$ geodesic, that
$h^{\star} - s^{\star} $ becomes precisely
$h^{\star}$.  Viewing this situation  as
in figure \ref{cayly}, we can now read off from figure
\ref{cayly}  that
\[ h^{\star}   - s^{\star} =
-\ln\left(\sech\left(\frac{a}{2}\right)\right) .
\]

Similarly notice, that 
 \[ -h^{\star}
+ t^{\star} =  \ln(\sech(l)), \]which as observed in figure \ref{quad} implies \[ - h^{\star} + t^{\star}  = \ln\left(\tanh\left(\frac{a}{2}
\right)\right) .\]

\begin{figure}[ht!]
\cl{\relabelbox\small\let\ss\scriptsize
\epsfxsize 2.2in \epsfbox{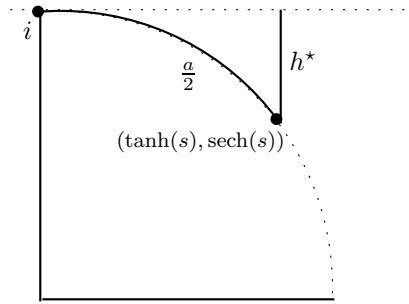}
\relabel {h}{$h^\star$}
\relabel {i}{$i$}
\relabel {a}{$a\over2$}
\relabela <-7pt,0pt> {t}{\ss$({\rm tanh}(s),{\rm sech}(s))$}
\endrelabelbox}  
\caption{Here we have placed the lower triangle
from figure \ref{quad}  into the upper-half plane model, sending the
ideal vertex to infinity and the $\frac{a}{2}$ segment on the
unit circle as pictured.   Recall that the  unit  circle in this
picture can be parameterized by hyperbolic distance from $i$
via $\tanh(s) +  \sech(s) i$.} 
\label{cayly} 
\end{figure}

With these computations we now have 
\[ a^{\star} = 2\left(
\ln\left(\tanh\left(\frac{a}{2}\right)\right) -
\ln\left(\sech\left(\frac{a}{2}\right)\right)\right)  =  2
\ln\left(\sinh\left(\frac{a}{2}\right)\right)\] 
as needed.
\endproof

\end{subsection}
\begin{subsection}{Convexity:  The proof of Lemma
\ref{cony}}\label{con} 

To prove $H$ is strictly concave, we start with the 
observation that  the objective function $H$  will certainly  
be a  strictly concave function on
$\frak{N}_x$ if the prism volume function $V(d^t(x))$ 
viewed as a function on   
\[ \{(A,B,C) \in (0,\pi)^3 \mid A + B + C < \pi \} \]  
turned out to be strictly concave. 
It is worth noting that, this implies  
$H$ is then strictly concave on all  of $\frak{N}$. 

There are several nice methods to explore the concavity of
$V(A,B,C)$. One could simply
check  directly that
$V$'s Hessian is negative  definite (as done in 
\cite{Le}), or one could exploit the visible injectivity of
the gradient, or one could bootstrap from the concavity 
of the ideal tetrahedron's volume. It is this last method
that will be presented here.
The  crucial observation is that any
family of ideal prism can be decomposed into three ideal
tetrahedra, as in figure
\ref{tet}.     So, we have 
\[ V(A,B,C) = \sum_{i=1}^3 T_i(A,B,C),\]
where $T_i$ is the volume of the $i^{th}$ tetrahedra in this
decomposition.
\begin{figure}[ht!]
\cl{\relabelbox\small
\epsfxsize 4in \epsfbox{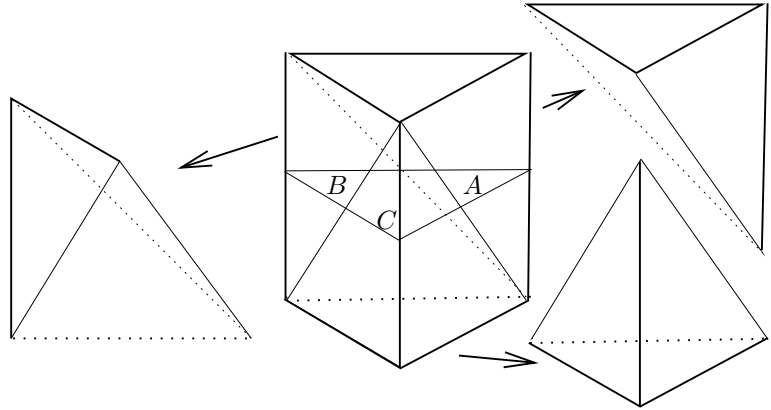}
\relabela <-3pt,0pt> {A}{$A$}
\relabela <-2pt,0pt>  {B}{$B$}
\relabela <-3pt,0pt> {C}{$C$}
\endrelabelbox}  
\caption{A decomposition of the ideal prism into
three ideal tetrahedra.   Notice that, the
angles in this decomposition are determined by the
affine conditions coming from the ideal vertices (see figure
\ref{horo1}) along with the condition that the angles meeting
along an edge slicing a prism face sum to $\pi$.  In particular, all the
angles depend affinely on the angles
$\{A,B,C\}$.  }
\label{tet} 
\end{figure}

Let us note some properties of the ideal tetrahedra
and its volume.  First, recall from figure \ref{horo1} that
 the dihedral angles corresponding 
to the edges meeting at a vertex of an
ideal tetrahedron are the angles of a Euclidean triangle.   In
particular, the fact that the  
constraint  $ \sum_{e \in v} \theta^e = \pi $ 
holds
at each vertex guarantees that 
an ideal tetrahedron is uniquely determined by any pair of
dihedral angles $\alpha$ and
$\beta$ corresponding to a pair of edges sharing a vertex. 
Furthermore, any pair of angles in
\[\{ (\alpha,\beta) \mid \alpha + \beta < \pi \}\] 
determines an ideal tetrahedron.
Note  the following fact (see \cite{Ri2}).

\begin{fact}
The ideal tetrahedron's volume function, $T(\alpha,\beta)$, is
strictly concave on the set 
\[\{ (\alpha,\beta) \mid \alpha + \beta < \pi \}\] 
and  continuous on this set's closure.  
\end{fact}

From figure \ref{tet}, each of the
$\alpha_i$ and
$\beta_i$ of the $i^{th}$ tetrahedron depend on the $(A,B,C)$
affinely. This fact immediately provides us with the
continuity assertion in Lemma \ref{cony}. To exploit the
tetrahedra's concavity  we will use the following straight-forward
lemma which follows immediately from the definition of concavity.

\begin{lemma}
Let $T$ be a strictly concave function
on the convex set $U \subset \Bbb{R}^{m}$ and for each
$i$ let  $L_i$ be an affine  mapping  from $\Bbb{R}^n$ to
$\Bbb{R}^{m}$ taking the convex set $V$ into $U$. Then
the function  
\[ V(\vec{x}) = \sum_{i = 1}^{k} T(L_i(\vec{x})) \]
is strictly concave on $V$  provided $L_1 \times \dots
\times L_{k}$ is injective.
\end{lemma}

Letting 
\[ L_i(A,B,C) = (\alpha_i(A,B,C),\beta_i(A,B,C))\]
we see that $V$ will satisfy the lemma if, for example, 
 the mapping 
\[(\alpha_1(A,B,C),\alpha_2(A,B,C)),\alpha_3(A,B,C))\]
is injective. Looking
at the decomposition in figure \ref{tet}, we see 
that we may in
fact choose
$\alpha_1(A,B,C) = A$,
$\alpha_2(A,B,C) = B$, and $\alpha_3(A,B,C) = C$. 
So, indeed,
we have our required injectivity and  $V(A,B,C)$ is strictly
concave as needed.

\end{subsection}
\begin{subsection}{Boundary control:  Proof of Lemma 
\ref{pushin}}\label{bound}

To begin proving Lemma \ref{pushin},  note that the  compactness of
$\frak{N}$'s
closure guarantees that 
$l(s)$ eventually hits the boundary again at some $y_1$ for some
unique   $s > 0$.   So we may change the speed of our line and assume  
we are using the line connecting  the two  boundary points, namely,\[l(s) = (1-s)y_0 +sy_1.\] 
In particular, lemma \ref{pushin}  is equivalent to knowing 
the following: for every pair of points $y_0$
and $y_1$ on $\partial \frak{N}$, with $l([0,1])$ unfoldable,  we have 
\[ \lim_{s \rightarrow 0^+} \frac{d}{ds}H(l(s)) \mbox{  } 
= \mbox{  } \infty. \]

Recalling that $H(d^t(x))  = \sum_{t \in \frak{T}} V(d^t(x))$, 
we see the lemma will follow if we demonstrate that
for any triangle 
 \[ - \infty < \lim_{s \rightarrow 0^+} \frac{d}{ds}V(d^t((s))) \leq \infty ,\]
 and for some triangle 
\[ \lim_{s \rightarrow 0^+} \frac{d}{ds}V(d^t((s))) = \infty . \]

The boundary is expressed in terms of angle data, so it would be nice to express
the $-2 \ln\left(\sinh\left(\frac{a}{2}\right)\right)$ coefficient in front
of the $dA^{\star}$ term in $dV$ (as computed in Section \ref{difff})
in terms of the angle data.  In fact, we can do even better and 
put this term in a form conveniently decoupling  the angle and curvature.

\begin{formula}\label{decup} 
$ - 2 \ln\left(\sinh\left(\frac{a}{2}\right)\right) $ is equal
to 
\[
\ln(\sin(B)) + \ln(\sin(C)) - \ln \left(\frac{ \cos(A - k^t(x)) - \cos(A)
}{k^t(x)}\right)  - \ln(-k^t(x)) . \]
\end{formula}
\proof[Proof of formula \ref{decup}]
This formula relies only on the hyperbolic law of cosines
which tells us 
\[ \cosh(a) = \frac{\cos(B) \cos(C) - \cos(A)}{\sin(A) \sin(B)}. \]
Using this relationship and the definition of $k^t(x)$  we now have
 \[ - 2 \ln\left(\sinh\left(\frac{a}{2}\right)\right)  = 
- \ln\left(\sinh^2\left(\frac{a}{2}\right)\right)   =-
\ln\left(\frac{\cosh(a) - 1}{2} \right) 
 \]   \[  = - \ln\left(\frac{\cos(B +C ) + \cos(A)}{\sin(B) \sin(C)} \right)
 = - \ln\left(\frac{- \cos(A - k^t(x)) + \cos(A)}{\sin(B) \sin(C)} \right),
\]  as needed.
\endproof

Using this formula, we will now enumerate the possible $y_0$ and the
behavior of $\frac{d}{ds} V(d^t((s)))$ in these various cases.  Let $F$
denote a finite constant.  We will be using the fact that if $L(s)$ is an
affine function of 
$s$ satisfying
$\lim_{s
\rightarrow 0^+}L(s) = 0$  then $\lim_{s
\rightarrow 0^+} \ln|\sin(L(l(s))|$ and the
 $ \lim_{s
\rightarrow 0^+} \ln|D(L(s))|$  can both be expressed as 
$\lim_{s \rightarrow 0^+} \ln(s) +F$.  Furthermore for convenience let 
$d^t(y_i)   = \{A_i,B_i,C_i\}$. 

\begin{enumerate}
\item
When   $d^t(y_0)$ contains no zeros and  $k^t(y_0) \neq
0$,  we have that $$\lim_{s \rightarrow 0^+} \frac{d}{ds} V(d^t((s)))$$ is
finite.

\item
When $d^t(y_0) = \{0,0,\pi\}$, we have 
 \[ \lim_{s \rightarrow 0^+} \frac{d}{ds} V(d^t((s))) = \]
\[ 
\frac{1}{2} \lim_{s \rightarrow 0^+}\ln(s)(\sigma^t(y_1 - y_0))- 
\frac{1}{2}\lim_{s
\rightarrow 0^+}\ln(s) (\sigma^t(y_1 - y_0)) +F  = F.\]

\item
In the case where   $d^t(y_0)$ contains zeros,  but  $k^t(y_0) \neq
0$, for each zero (assumed to be  $A_0$ below) 
 we produce a term
in the form
 \[ \lim_{s \rightarrow 0^+} \frac{d}{ds} V(d^t((s))) = 
\lim_{s \rightarrow 0^+} \left( -\ln(s)(A_1 - A_0) \right),\]
plus some finite quantity.

\item
When $k^t(y_0) = 0$  and no angle is zero
 \[ \lim_{s \rightarrow 0^+} \frac{d}{ds} V(d^t((s))) = 
\lim_{s \rightarrow 0^+} \frac{1}{2}\ln(s)
(\sigma^t(y_1 - y_0)) +F .\]

\item
When $k^t(y_0) = 0$  and one angle, say $A_0$, in $d^t(y_0)$ is zero we have 
 \[ \lim_{s \rightarrow 0^+} \frac{d}{ds} V(d^t((s))) = 
-2\lim_{s \rightarrow 0^+}\ln(s)(A_1 - A_0)) + \lim_{s \rightarrow 0^+}
\ln(s) (\sigma^t(y_1 - y_0)) +F .\]

\end{enumerate}

So, the first two cases produce finite limits.  In order to  
understand the next three  limits we make some simple observations.  First, 
if  $A_0 =0$ and $l(s) \bigcap \frak{N} \neq \phi$, then $A_1-A_0 > 0$.  So,
limits from the third case evaluate to $+ \infty$.     Secondly, note
that   when 
$k^t(y_0) = 0$ and $l(s) \bigcap \frak{N} \neq \phi$,  that $ \sigma^t(y_1 -
y_0) =A_1+B_1+C_1 -( A_0+B_0+C_0 ) < 0$ and, hence,
 the limits from the  fourth
case are $+\infty$ too.  Combining these observations we see the
fifth case always produces a$+\infty$ limit as well.

For each triangle, the answer is finite or positive infinity.  
So, all we need to do is guarantee that for some triangle we achieve
$+\infty$. To do this, note that in order for $y_0$ to be unfoldable,
there is some triangle $t$ such that   
$d^t(y_0) = \{A_0,B_0,C_0\} \neq \{0,0,\pi\}$, however, either $k^t(y_0) =0$ or
some angle is zero.  So we have at least one triangle in case 3,4, or 5,
 as needed.

\end{subsection}
\end{section}

\begin{section}{Feasibility: 
The proof of Proposition \ref{emf} }\label{sec:4}

We will first demonstrate that,
 if $\frak{N}^p$ is non-empty for a
Delaunay 
$p \in \Bbb{R}^{|E|}$, then  $p$ is feasible. 
In order to accomplish this, we first derive a useful formula.
To state the formula, we let  $O(S)$ denote the outside of $S$; in other words, let $O(S)$ be the
set of all pairs $(e,t)$ such that
$e$ is an edge of exactly one face in $S$, and $t$ is 
the face in $F-S$ which contains $e$.

\begin{formula}\label{sett}
For any   set of triangles $S$ \[  \sum_{e \in S}  \theta^e =
\sum_{t
\in S} \left( \pi -
\frac{k^t}{2}\right)  + \sum_{O(S)}
\left(\frac{\pi}{2} - \psi^e_{t}\right) . \] 
\end{formula}
\proof
\[ \sum_{E } \theta^e = 
\sum_{e \in S} (\pi - \psi^e) 
= \sum_{e \in S} 
\left( \left(\frac{\pi}{2} - \psi_{t_1}^e \right) 
+ \left(\frac{\pi}{2} - \psi_{t_2}^e \right)\right) \]
\[ = \sum_{t \in S} 
\left( \pi  + \frac{ \pi - \sigma^t }{2}\right) 
+ \sum_{O(S)}  
\left(\frac{\pi}{2} - \psi^e_{t} \right) 
. \]
  Substituting the definition of 
$k^t$ gives the needed formula.
\endproof

Since $\frak{N}^p$ is non-empty, $\frak{N}^p$ contains a Delaunay angle system
$x$.
Apply the right and left hand sides of formula \ref{sett} to this $x$.
Since  $- k^t(x) > 0$ and 
 $\frac{\pi}{2} - \psi^e_{t}(x) >0$ (as in the proof of
  lemma \ref{delun}),   
removing these terms from the right hand side 
strictly reduces the right hand side's size and implies   that 
\[   \sum_{e \in S}  \theta^e(x) > \pi |S|, \]
as needed.
  
To demonstrate the converse, 
that a feasible Delaunay $p \in \Bbb{R}^{|E|}$
has a non-empty 
 $\frak{N}^p$,  we 
 will characterize the points in $\frak{N}^p$ as the 
 minimal elements in the larger set  
\[\frak{N}^{+}_p=\{x \in \Bbb R^{3|F|}
\mid \alpha^e_t(x) \in (0,\pi) , \theta^e_t(x) \in (0,\pi),\Psi(x)=p, \,\forall (e,t)\}\]
with respect to the  objective function 
\[m(x)=
\left( \sum_{\{t \mid k_t(x) > 0 \}} k_t(x),
 |\{t :  k_t(x)=0\}|,|\{(e,t) : \alpha_e^t(x)=0\}| \right) \]
which takes values in 
$\Bbb R \times \Bbb Z \times \Bbb Z$, given its dictionary order.
By construction, $y \in \overline{\frak{N}^{+}_p} $ satisfies
$m(y)=(0,0,0)$  if and only if $y \in \frak{N}_p$.

The key fact is that $\frak{N}^{+}_p$
is automatically non-empty, since $p$ being Delaunay allows us to
determine a point $x \in \frak{N}^{+}_p$ by setting
  $\theta^e_t(x) = \frac{\theta^e}{2}$. By compactness of
$\overline{\frak{N}^{+}_p}$ there is a subset of 
$\overline{\frak{N}^{+}_p}$  upon  which $m$ achieves its minimal 
value.  Hence, the following lemma, proved in Section \ref{surj}, 
will finish our proof of Proposition \ref{emf}.

\begin{lemma}\label{linmin}
The minima of $m$ on  $\overline{\frak{N}^{+}_p}$ is $(0,0,0)$.
\end{lemma}

\begin{subsection}{Minimality: The proof of Lemma \ref{linmin}}\label{surj}

The proof of Lemma \ref{linmin} involves the careful exploration of
several ``non-local"  objects.  This exploration, however, 
relies on  some  simple local observations, collected in the following sub-lemma.

\begin{sublemma}\label{tech}
Suppose $x \in \overline{\frak{N}^+_p}$ for a feasible Delaunay $p$.
\begin{enumerate}

\item \label{s1}
 $k^t(x)<0$,
 then for every $e \in t$ we have $\theta^e_t(x)>0$.
\item \label{s2}
Let the edges of $t$ be $\{e,e_1,e_2\}$, $k_t(x) \geq 0$,
 and  $\alpha^e_t(x)=0$, then
\begin{enumerate}
\item[\rm(a)] 
$\alpha^{e_i}_t(x) >0$.
\item[\rm(b)]  
$\theta^{e_i}_t(x)>0$.
\end{enumerate}
\item \label{s3}
It is impossible for pair of triangles $t_1$ and $t_2$ 
satisfying $k_{t_i}(x) \geq 0$ and  sharing the edge $e$ to satisfy
$0 = \alpha^e_{t_1}(x)=\alpha^e_{t_2}(x)$
\item  \label{s0}
\[ \frak{N}^+_p = \{x \in \Bbb R^{|E|} \mid \alpha^e_t(x)>0,\theta^e_t(x)>0\}\]
\end{enumerate}
\end{sublemma}
\proof
Part \ref{s1} follows from the following fact that
$\theta^e_t = \alpha^e_t-k^e_t/2$.

To begin part \ref{s2}, it is useful to note that 
the above formula for $\theta^e_t$ demonstrates that 
$k^t(x)=0$, since $\theta^e_t(x) \geq 0$.  
Part 2a now follows from the fact that, if there are two zeros 
and $k^t(x) = 0$, then the edge $e$ where the two zeros
live satisfies
\[ \theta^e(x)  = \theta^e_t(x) +\theta^e_s(x)
  = \pi+\theta^e_s(x)  \geq \pi ,\]
contradicting the Delaunay constraint.
Similarly 2b follows  by noting that 
if $\theta^{e_1}_t(x) = 0$,  then 
\[ \theta^{e_1}_t(x) + \theta^e_t(x) =\pi- \alpha^{e_2}_t(x)  = 0 , \]
hence, 
$\alpha^{e_2}_t(x)  = \pi$ and $\theta^{e_2}(x) \geq \pi$.

Part \ref{s3} follows by noting that, under these circumstances
 $\theta^e_{t_1}(x) = 0$ and $\theta^e_{t_2}(x) = 0$, 
hence, $\theta^e(x) = 0$, contradicting the Delaunay constraint. 

Part \ref{s0} follows from the
fact that, $\pi-\alpha^e_t =\theta^{e_1}_t+\theta^{e_2}_t$, where 
$\{e,e_1,e_2\}$ index the three edges of $t$, together with the fact that 
$\theta^{e}_{t_1}(x)+\theta^{e_2}_{t_2}(x) = \theta^e(x)<\pi$.\endproof

We will now prove a weak version of Lemma \ref{linmin}.

\begin{lemma}\label{vacuum}
The minima of $m$ on  $\overline{\frak{N}^{+}_p} $ is in the from $(0,0,z)$.
\end{lemma}
\proof
Let $x$ be  a minimal point of $m$ in $\overline{\frak{N}^{+}_p}$.
Let $Z$ denote the  set of 
all triangles with $k^t(x) \geq 0$.
 This sub-lemma  is equivalent to $Z$  being empty.
We  will hypothesize that   $Z$ is not empty and produce a contradiction.

The first observation needed about 
$Z$ is that $\tot{Z}$ is not all of $M$.
This follows from the feasibility condition  on $p$, namely,
\[\sum_{t \in F} k^t(x) = 2(\pi |F| - \sum_{e \in S} \theta^e(x)) < 0, \]
hence,  there is negative curvature somewhere.

By the minimality   of $x$,  there can be no conformal transformation
capable of  moving negative  curvature into  $Z$. 
In particular, for any 
$\frak{c} \in \frak{C}$ that satisfies  
$m(x+\frak{c}) < m(x)$'s, we have
$x+ \frak{c}$ is in the complement of  $\overline{\frak{N}^{+}_p}$. 
 Hence, by sub-lemma \ref{tech}  part \ref{s0}, we have that either  
$\theta^e_t(x+\frak{c})<0$ or $\alpha^e_t(x+\frak{c})<0$ for some $(e,t)$.
To exploit this fact
we will attempt to conformally vacuum some curvature out of $Z$. 
Suppose we are at  an internal
boundary edge, $e_0$, of $Z$, ie $(e_0,t_{-1}) \in O(Z)$. Denote by
$t_0$  the triangle in $Z$  containing $e_0$.
First, we will attempt to move negative  curvature our of $Z$ 
with a conformal
transformation in the form $\epsilon w_{e_0}$.
Since  $t_{-1}$ has  negative  curvature, 
by sub-lemma \ref{tech}  part \ref{s1}, 
the obstruction 
being able to do so
must be due to a  zero angle along $e_0 \in t_0$.  
Furthermore, since $k_{t_0}(x) \geq 0$,
 by  sub-lemma \ref{tech}  part 2a,  $t_0$ contains only one 0.Now
we will continue the attempt to remove curvature. Let $e_1$ be the other
edge sharing the unique zero angle along 
$e_0$ in $t_0$; if  $e_1$ is another boundary of $Z$ edge we stop. If
$e_1$ is not a boundary edge of $Z$ use the transformation
$\epsilon(w_{e_1} + w_{e_0})$ 
to continue the effort
to  remove curvature.  
 By sub-lemma \ref{tech}  part 2b, 
we are in the same position as above and there is 
a unique zero on  the non $t_0$ side of $e_1$, call this face $t_1$. 
Notice,  we can continue  this  process and form a uniquely determined set of
 edges  and triangles, $\{e_i\}$ and $\{t_i\}$,   and  a sequence of conformal  
transformations $\epsilon \sum_{i =0}^n w_{e_i} \in \frak{C}$,
 see figure \ref{vac}.
Since such a sequence can never  intersect itself
with out producing a pair triangles as in the hypothesis of sub-lemma
 \ref{tech} part \ref{s3}, the sequence  must eventually poke 
through $Z$ into $Z^c$.
We will refer to the set of faces arising in this construction as a 
{\it vacuum}, see figure \ref{vac}.

\begin{figure}[ht!]
\cl{\relabelbox\small\let\ss\scriptsize
\epsfxsize 4.2in \epsfbox{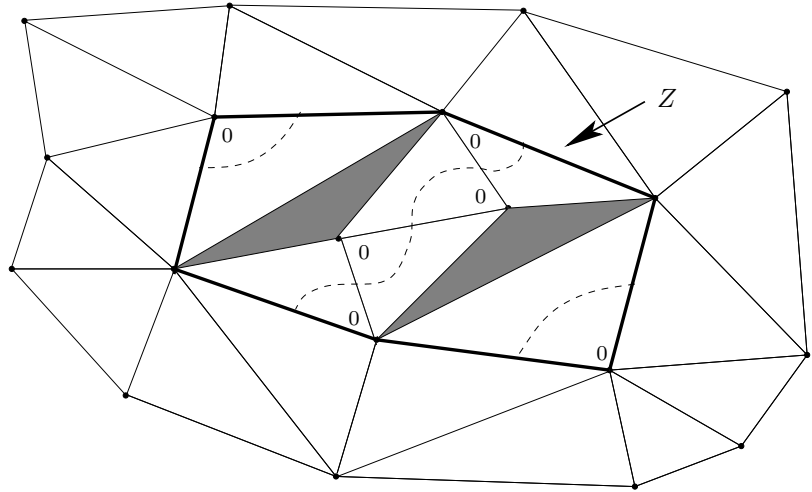}
\relabel {0}{\ss$0$}
\relabela <0pt,-2pt>  {01}{\ss$0$}
\relabel {02}{\ss$0$}
\relabel {03}{\ss$0$}
\relabel {04}{\ss$0$}
\relabel {05}{\ss$0$}
\relabela <-2pt,2pt> {Z}{$Z$}
\endrelabelbox}  
\caption{In this figure the region of zero curvature, $\tot{Z}$,  
is the set with the bold boundary as labeled.   We have also attached to $Z$
all its vacuums and indicated them with the squiggly lines.  As we can see in 
the figure, these vacuum transformations fail  to remove curvature from
$Z$, just  as in the proof of lemma \ref{vacuum}.  Notice, the set  $S$,
$Z$ minus its vacuums,  is in this case the pair of shaded triangles. }
\label{vac}
\end{figure} 

Now in order to utilize our vacuums, 
we need to remove them.
Let  $S$ be the removal from $Z$ of all these vacuums, see figure \ref{vac}.   
The key is the following.

\begin{sublemma}\label{vacin}
If $Z$ is non-empty, then $S$ is non-empty.
\end{sublemma}
\proof
First, note two distinct vacuums can never  share an edge.   
To see this, call a  vacuums {\it side boundary} any edge  of a face in 
a vacuum  facing a zero.   Now simply note, if the intersection of  two vacuums 
contains an edge, then it contains a first edge $e_i$ with respect to one of 
the vacuums. There are two possibilities for this edge. One is that 
$t_{i+1}$ has two zeros  and $k^t(x) \geq 0$, contradicting sub-lemma
 \ref{tech}  part 2a, 
 or that $e_i$ is a side boundary 
of both  vacuums. In this case, we produce a pair of triangles satisfying 
the hypothesis of sub-lemma \ref{tech}  part \ref{s3},
which also can not exist.

Since vacuums cannot intersect themselves or  share edges with
distinct  vacuums,  $S$  would be non-empty if the side
boundary of a vacuum 
had to be  in $Z$'s interior.  
 Look at any side boundary edge $e$ of a fixed vacuum.
Note $e$  cannot be in $Z$'s boundary because, in  this case, the
vacuum  triangle it belonged to would have at least two zeros 
 and $k^t(x) \geq 0$,
contradicting sub-lemma \ref{tech}  part 2a.  
Hence, $S$ is nonempty.
\endproof

With our non-empty $S$ in hand, we may 
now finish lemma \ref{vacuum} by looking  
at formula \ref{sett} with respect to $S$.
Every edge in the boundary of $S$ faces a zero on its $S^c$
side
in a triangle with $k^t(x) \geq 0$ (see figure \ref{vac} once again), so
\[ \sum_{O(S)}   
\left(\frac{\pi}{2} - \psi^e_{t}(x)\right) \leq 0.\]
Also, each triangle in $S$ has non-negative
curvature, hence, from   formula \ref{sett},  
we have violated  the feasibility condition, as needed.
\endproof

We can now bootstrap from  Lemma \ref{vacuum} to Lemma
\ref{linmin}. Once again, assume
 $x$ is  a minimal point of $m$ in $\overline{\frak{N}^{+}_p} $.
From Lemma \ref{vacuum},  we know that $m(x)=(0,0,z)$.
We will suppose this  minimal $z$ is strictly greater 
than zero and produce a contradiction.

Since $z >0$ there is an  angle slot which is zero.
Fix such an  angle slot, $(e,{t_0})$, and let $\{e,e_0,e_1\}$
index the edges of $t_0$.
Just as in lemma \ref{vacuum}, we now may use 
sub-lemma \ref{tech}  part 2a to prove: in order for $x$ 
not to be conformally equivalent to a point 
where $m$ is decreased  by the $\epsilon w_{e_1}$ transformation, 
 we must have a zero along the 
$e_1$ edge of $t_0$'s neighbor $t_1$.
Call  this neighboring triangle $t_1$ and let  $e_2$ be
another edge in $t_1$ next to a  zero angle slot in $t_1$.
Since 
$k^{t_1} < 0$, we may repeat the above procedure, letting $e_2$ 
play the role of
 $e_1$  and constructing an 
$e_3$ and $t_2$.
We may continue this process  indefinitely
forming a sequence of edges  $\{e_i\}_{i=1}^{\infty}$.
Since there are a finite number of edges, eventually
the  sequence $\{e_i\}_1^{\infty}$  will have some $k<l$ such that 
$e_k = e_l$  and
$e_{k+1} =e_{l+1}$,  see figure \ref{sna}. 
I'll call such a finite set of edges  an {\it accordion},
see in figure \ref{sna}  for an example. 

\begin{figure}[ht!]
\cl{\epsfxsize 4.2in \epsfbox{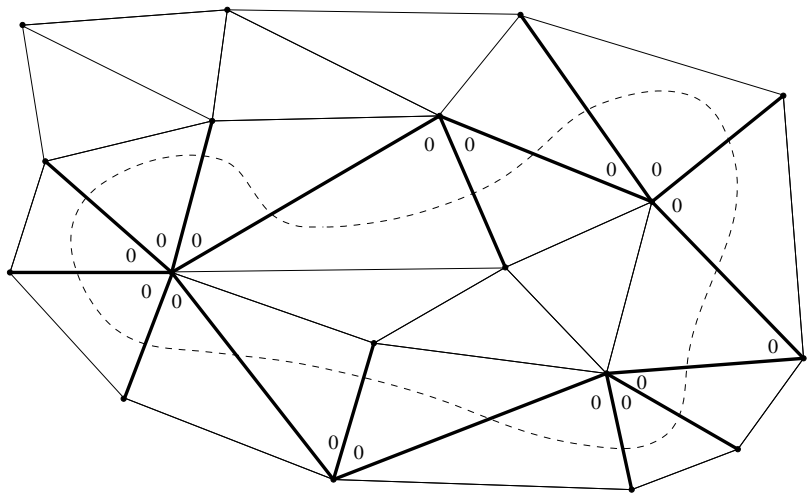}}  
\caption{Here we have an accordion,
the set of edges forming a loop indicated  with the squiggly
line. Notice, to geometrically realize the 
 zero angles in the pictures corresponds to squeezing the 
accordion.  The algebraic inability to squeeze various accordions is
at the heart of the proof of Lemma \ref{linmin}.}  
\label{sna} 
\end{figure}

Now we can produce a contradiction to this accordion's existence. 
If we denote the three edges of  $t_i$ as $\{e_i,e_{i+1},e_{i,i+1}\}$, 
then we have that $\alpha^{e_{i,i+1}}_{t_i}(x)=0$ for each $i$. 
This, together with the fact that 
$\alpha^{e_{i,i+1}}_{t_i} = \psi_{t_i}^{e_i} + \psi_{t_i}^{e_{i+1}}$,
demonstrates that 
\[ 0 = \sum_{i = k}^{l-1} \alpha^{e_{i,i+1}}_{t_i}(x)
  =  \sum_{i=k}^{l-1} \psi^{e_i}(x). \]
However, by the Delaunay hypothesis on $x$, we have that  $\sum_{i=k}^{l-1} \psi^{e_i}(x) > 0$, our needed contradiction.

\end{subsection}

\end{section}

\begin{section}{Generalizations, consequences, and questions}\label{sec:5}

\begin{subsection}{Straightforward generalizations}

The techniques in the previous sections apply nearly immediately to certain 
simple generalizations.  The most important of these  is to generalize  the 
use of a  triangulation to keep track of the informal angle complements.  
 The following structure will prove much more useful.
A  polygonal CW--complex will be called a {\dfn{polygonal decomposition}}
 if it  satisfies the following 
property: 
we may triangulate each polygon using 
its vertices in such way that  the resulting CW--complex lifts to a
triangulation
of the
surface's universal cover.
A geodesic  polygonal decomposition  will be
called a
{\dfn{Delaunay decomposition}}  
if it can  be 
constructed  using Delaunay's
empty sphere method, as described in
the
Delaunay Construction.

After replacing the triangulation in Section
\ref{sec:2} with a 
polygonal decomposition, we have the following
scholium to Theorem \ref{cir1}.
The details of how to accomplish this can
be found in \cite{Le}.

\begin{scholium}\label{cir2}
Suppose our
surface is compact and boundaryless.  
A non-sing\-ular $p \in \Bbb R^{|E|}$
is realized by a unique Delaunay decomposition if and only if
$p$ is
feasible and Delaunay.
\end{scholium}

An even more straight
forward generalization is to to allow singularities and
 boundary with
corners.  
Since the argument is localized to triangles, we may freely rid
the 
non-singularity condition with the price that our hyperbolic surface
may have specified singular points.
The case of a surface with boundary
follows by doubling the surface.
More precisely,  at a
 boundary edge $e
\in t$ of an
geodesic triangulation, we may let $\psi^e$ be defined  as
 $\psi^e_t$.  If we modify the Delaunay
condition to include the hypothesis that $\psi^e \in
(0,\frac{\pi}{2})$ for every boundary edge $e$, then 
Theorem \ref{cir1} still 
holds when  we remove both the non-singularity hypothesis and  
the implicit boundaryless condition.

\end{subsection}

\begin{subsection}{A consequence}

We are now in a position to discuss Theorem \ref{cir1}'s
implications about the Teichm\"uller space of a surface with
 genus greater
than one and  at least one distinguished point. 
First, let us state the view of  Teichm\"uller  space  that we will be
exploring.  Fix a differentiable surface $N$ and a
specified ordered set of points $\{v_i\}_{i=1}^{|V|}$ on $N$.
Let $(G,V,f)$ denote a triple consisting of a  closed hyperbolic
surface $G$, an ordered set of vertices $V$ in $G$, and an orientation
preserving diffeomorphism
$f$ from $N$ to $G$ such
that $f\left(\{v_i\}_{i=1}^{|V|}\right)=V$ as an ordered set.
As a set, let the Teichm\"uller space be
the equivalence classes  of triples $(G,V,f)$, where
 $(G,V,f)$
and $(\overline{G},\overline{V},\overline{f})$ will be considered
equivalent
if there is
an isometry $i$ from $G$ to $\overline{G}$ such that   $i$ is isotopic to
$\overline{f} \circ f^{-1}$. 
We will now use Scholium \ref{cir2} to form a natural tessellation of
Teichm\"uller space.  

\medskip
{\bf Context}\qua
A combinatorially equivalent version of the  tessellation of the
Riemann Teichm\"uller space that
we will construct  was originally originally introduced in
\cite{Ha}.  There have been
several beautiful  methods for geometrically realizing this tessellation
related to our method.   One of these methods starts with the wonderful observation that the  Riemann Teichm\"uller space   can be related to the
Teichm\"uller space of singular Euclidean structures, which was first
accomplished in \cite{Tr}.  With this view, one can then exploit the fact
that
the Teichm\"uller spaces of singular Euclidean structures has a
particularly
natural cellular structure, as in \cite{Bo} or \cite{Ve}.  In fact, this
was the
setting where   the Euclidean version of  Theorem \ref{cir1} was
initially put to use,  see  Bowditch's \cite{Bo}.  To geometrically
accomplish this
identification,  one first identifies the Riemann  Teichm\"uller space with
the
space of  complete finite area hyperbolic structures with the
distinguished
points viewed as cusps. Using some form of horoball decoration (as in
\cite{Pe}) one can relate this  hyperbolic structure directly to a
singular
Euclidean structure using the construction  in \cite{Ep}.  Here we will
take a
different geometric view in order to construct the tessellation.  We
will view   the Riemann Teichm\"uller space as the space of hyperbolic
structures
with distinguished points.  We shall find that  the Delaunay
decomposition
and  Scholium \ref{cir2}  provide us with the natural  polyhedral  cells
of the Teichm\"uller space's tessellation. Due to the large number of
excellent references, the discussion of the construction of the tessellation 
below will be kept relatively terse.

\medskip{\bf Tessellation construction}\qua
First, we will show that the Teichm\"uller space can be
identified, as a set, with a union of the polyhedra described in Scholium
\ref{cir2}.  We may
restrict our attention to the polygonal decompositions associated to
non-singular data.   Let $\bf{P}$ be the set of such  polygonal
decompositions of $N$ viewed up to isotopies fixing
$\{v_i\}_{i=1}^{|V|}$. To be explicit, the elements of $\bf{P}$ will be
equivalence classes of pairs  $(P,f)$ where
$P\in \bf{P}$ with an order on its
vertices, $V$, and $f$ is a homeomorphism of $N$ to  $\tot{P}$ such that
$f\left(\{v_i\}_{i=1}^{|V|}\right)=V$  as ordered sets.  Let  $(P,f)$ and
$(\bar{P},\bar{f})$ be equivalent if there is an
orientation preserving homeomorphism $I$ of $N$ which is isotopic to the
identity relative 
$\{v_i\}_{i=1}^{|V|}$ and a CW equivalence $h$ of $P$ and $\bar{P}$, 
such that $I \circ f^{-1}  =
\bar{f}^{-1}  \circ h$ when restricted to every cell.   
 For each element $[(P,f)]$ of
$\bf{P}$, fix an order on its edges and let  $E([(P,f)])$ be  the open
convex polyhedron determined by the
non-singular  points   in $(0,\pi)^{|E|} $ satisfying  the
feasibility  condition. 
Notice,  each
point in Teichm\"uller space will correspond to at least one  point in
$E(\bf{P}) = \bigcup_{[(P,f)] \in
{\bf P}} E([(P,f)])$ via the Delaunay triangulation.  In fact, this
correspondence is a well defined mapping, since among the
decompositions in $\bf{P}$ there is no
way to isotopically interchange the edges without interchanging the
vertices as well. By the  uniqueness part of Scholium 
\ref{cir2} this mapping is injective, and the existence assertion in
Scholium \ref{cir2} guarantees that this mapping  is surjective.  So,
as a set, we may identify the
Teichm\"uller space  with $ E(\bf{P})$.

Notice that,  whenever $\theta^e$ degenerates
to a  $\pi$, that the edge $e$ in the underlying polygonal decomposition is
removed in such a way that, after its removal,  we arrive at a
new polygonal decomposition with all the other data agreeing.  We may use
this observation  to glue up  $\bigcup_{[(P,f)] \in {\bf P}}
E([(P,f)])$ in way that preserves the Euclidean structure
 of each $E([(P,f)])$ piece   and, hence, gives $E(\bf{P})$  the structure
of
an open Euclidean cone manifold.  From our set correspondence 
in the above construction,  we have
provided the Teichm\"uller space with a
Euclidean cone manifold structure, and, in particular, the structure of a
$E - V = 3\chi(N) + 2V=6g-6+2|V|$ dimensional manifold.  

Now, by looking at
the conformal structure of $G-V$, we have a bijective correspondence
between the points in  $E(\bf{P})$ and the  Riemann Teichm\"uller  spaceof a surface
with $|V|$ punctures, and we shall call this mapping $F$.  As a scholium to
the proof of Proposition
\ref{tri1}, as our
coordinates in $E(\bf{P})$ continuously vary, we continuously
vary the shapes of the 
geodesic polygons which are being used to build our hyperbolic
 surfaces, and in
particular, are quasi-conformally
deforming the conformal structure on the corresponding surfaces.  
So $F$ is a
bijective continuous map between
$6g-6+2|V|$ dimensional manifolds, hence, a homeomorphism due to the
invariance of domain.
Notice that, the mapping
class group respects  our   $E([(P,f)])$ cells, and so we have
formed a natural tessellation of Teichm\"uller space which descends to a
natural  tessellation of Moduli space.  We sum up this discussion with the
following result.

\begin{corollary}\label{tess}
The Euclidean cone manifold  $E(\bf{P})$ is homeomorphic to the
Teichm\"uller space of a surface with genus greater
than one and $|V|$ distinguished points.  Furthermore, the  Euclidean 
$E([(P,f)])$ cells provide a tessellation of Teichm\"uller space
invariant under the action of the mapping class group.
\end{corollary}

\end{subsection} 

\begin{subsection}{Questions}

The use of hyperbolic volume in this paper to solve Proposition \ref{tri1} 
and, hence,
Theorem \ref{cir1}, could conceivably be used to prove  the analogous
questions
directly in the spherical case. The prism, constructed in Section \ref{Hyp}, becomes 
the 
following twisted
prism, also observed as the right object for this game by
Peter Doyle. 

\medskip
{\bf Twisted prism construction}\qua Fix a
point $p$ in $H^3$ viewed as the origin in the Poincare
ball model  of $H^3$.  Imagine
 $S^{\infty}$ is given the unit sphere's induced metric from $\Bbb{R}^3$.
To construct the {\dfn{twisted prism}},
 place  a geodesic triangle on the sphere
and take the convex hull of this triangle's vertices together with
the point $p$. 

 To construct the associated objective function, once again,
sum up the volumes of the twisted prisms associated to an
angle system with  positive curvature.   It
is easy to see that the critical points of this
objective function are once again
precisely
the uniform angle systems, but
the objective function fails to be convex. Can this objective function
still be used
to arrive at the analogs of Proposition \ref{tri1} and Theorem
\ref{cir1}?

A far reaching generalization of the Thurston--Andreev Theorem, due to He (see
\cite{He2}), completely characterizes finite sided convex hyperbolic  polyhedra.
Along with the techniques in \cite{He1}, this generalization could be used to characterize generalized versions of the convex polyhedra constructed here.
The hyperbolic volume technique utilized here can be set up in this more general 
setting, though, as in the spherical case, a certain  amount of control is lost.
Can the hyperbolic volume techniques explored in this paper
be used to prove these more general results?

\end{subsection}

\end{section}

\end{document}